\documentclass[english,12pt,titlepage,bibliography=totoc,headsepline]{article}

\usepackage[arrow, matrix, curve]{xy}
\usepackage{graphicx}
\usepackage{amsmath} 
\usepackage{amssymb} 
\usepackage{amsfonts}
\usepackage{empheq}
\usepackage{amsthm}
\usepackage{makeidx} 
\usepackage{verbatim}
\usepackage{epsfig}
\usepackage{hyperref}
\usepackage[T1]{fontenc}
\usepackage[latin1]{inputenc} 
\usepackage{bbm}

\usepackage{lmodern}		
\usepackage{textcomp}

\numberwithin{equation}{section}

\newcommand{\R} {{\mathbb R}}
 
\newcommand{\N} {{\mathbb N}}
\newcommand{\Z} {{\mathbb Z}}

\newcommand{\grad}{\operatorname{\grad}} 

\DeclareMathOperator{\Lip}{Lip}
\DeclareMathOperator{\aco}{aco}
\DeclareMathOperator{\diam}{diam}
\DeclareMathOperator{\image}{Im}

\DeclareMathOperator{\rank}{rank}
\DeclareMathOperator{\codim}{codim}

\renewcommand{\epsilon}{\varepsilon} 
\renewcommand{\rho}{\varrho} 
\renewcommand{\phi}{\varphi} 
\newcommand{\e}{\varepsilon} 
 
\renewcommand{\b}{\beta} 
\renewcommand{\a}{\alpha}

\newcommand{\sgn}{{\rm sign }}

\theoremstyle{plain}  
\newtheorem{prop}{Proposition}[section]  
\newtheorem{coron}[prop]{Corollary}
\newtheorem{lemma}[prop]{Lemma} 
\newtheorem{theo}[prop]{Theorem}

\theoremstyle{definition}
\newtheorem{bsp}{Example}
\newtheorem{rem}{Remark}  

\renewcommand{\proof}{\noindent \textit{Proof.} }  

\title{Entropy numbers of convex hulls in Banach spaces and applications} 

\author{Bernd Carl \\
Mathematisches Institut, Universit\"at Jena\\
Ernst-Abbe-Platz 2, 07743 Jena, Germany\\
email:\ bernd.carl@uni-jena.de\\
\qquad
\\
Aicke Hinrichs \\
Mathematisches Institut, Universit\"at Jena\\
Ernst-Abbe-Platz 2, 07743 Jena, Germany\\
email:\ a.hinrichs@uni-jena.de\\
\qquad
\\
Philipp Rudolph \\
Mathematisches Institut, Universit\"at Jena\\
Ernst-Abbe-Platz 2, 07743 Jena, Germany\\
email:\ philipp.rudolph@uni-jena.de}

\begin{document}

\maketitle

\begin{abstract} 

In recent time much attention has been devoted to the study of entropy of convex hulls in Hilbert and Banach spaces and their applications in different branches of mathematics. In this paper we show how the rate of decay of the dyadic entropy numbers of a precompact set $A$ of a Banach space $X$ of type $p$, $1 < p \le 2$,  reflects the rate of decay of the dyadic entropy numbers of the absolutely convex hull $\aco(A)$ of $A$. Our paper is a continuation of the paper \cite{CHP11}, where this problem has been studied in the Hilbert space case. We establish optimal estimates of the dyadic entropy numbers of $\aco(A)$ in the non-critical cases where the covering numbers $N(A,\e)$ of $A$ by $\e$-balls of $X$ satisfy the Lorentz condition
$$ \int_0^{\infty} \left( \log_2(N(A,\e)) \right)^{s/r} \, \mathrm{d}\e^s < \infty $$
for $0<r<p'$, $0<s<\infty$ or 
$$ \int_0^{\infty} \left( \log_2(2+\log_2(N(A,\e))) \right)^{-\a s} \, \left( \log_2(N(A,\e)) \right)^{s/r} \, \mathrm{d}\e^s < \infty $$
for $p'<r<\infty$, $0<s \le \infty$ and $\a \in \R$, with the usual modifications in the case $s=\infty$. The integral here is an improper Stieltjes integral and $p'$ is given by the Hölder condition $1/p + 1/p' = 1$. It turns out that, for fixed $s$, the entropy of the absolutely convex hull drastically changes if the parameter $r$ crosses the point $r = p'$. It is still an open problem what happens if $r=p'$ and $0<s<\infty$. However, in the case $s=\infty$ we consider also the critical case $r=p'$ and, especially, the Hilbert space case $r=2$. 
   
We use the results for estimating entropy and Kolmogorov numbers of diverse operators acting from a Banach space whose dual space is of type $p$ or, especially, from a Hilbert space into a $C(M)$ space. In particular, we get entropy estimates of operators factoring through a diagonal operator and of abstract integral operators as well as of weakly singular convolution operators. Moreover, estimates of entropy and Kolmogorov numbers of the classical and generalized Riemann-Liouville operator are established, complementing and extending results in the literature. \medskip

 \medskip


\end{abstract}

\section{Introduction} 

Let $(M,d)$ be a metric space and $B_M(x,\e) := \left\{ y \in M \mid d(x,y) \le \e \right\}$ the closed ball with center $x$ and radius $\e$. For a bounded subset $A \subset M$ of $M$ and a natural number $n \in \N$, the \textit{$n$-th entropy number of $A$} is defined by 
$$ \e_n(A) := \inf \left\{ \e \ge 0 \, \middle| \, \exists \, x_1, \ldots, x_n \in M:  A \subset \bigcup_{i=1}^n B_M(x_i,\e) \right\}. $$
Moreover, the \textit{$n$-th dyadic entropy number} of $A$ is given by $ e_n(A) := \e_{2^{n-1}}(A)$. The entropy numbers permit us to quantify precompactness, the rate of decay of $\e_n(A)$ can be interpreted as a degree of precompactness of the set $A$. In order to apply this idea to linear bounded operators $T \in \mathcal L(X,Y)$ between Banach spaces $X$ and $Y$ we define the \textit{$n$-th entropy number of $T$} by 
$$ \e_n(T:X \to Y) := \e_n(T(B_X)), $$ 
where $B_X$ is the closed unit ball of $X$. The $n$-th dyadic entropy number of $T$ is given by $e_n(T):= \e_{2^{n-1}}(T)$. According to the definition, the operator $T$ is compact if and only if $\e_n(T) \to 0$ for $n \to \infty$. The speed of decay of the entropy numbers of $T$ can be seen as a measure for the compactness of $T$. The concept of \textit{covering numbers} is closely related to that of entropy numbers. For a bounded subset $A \subset M$ of a metric space $M$ and $\e > 0$ we define the $\e-$covering number $N(A,\e)$ by 
$$ N(A,\e) := \min \left\{ n \in \N \, \middle| \, \exists \, x_1, \ldots, x_n \in M:  A \subset \bigcup_{i=1}^n B_M(x_i,\e) \right\}. $$

In order to formulate precise estimates, we present a generalization of \textit{Lorentz sequence spaces}. These spaces will play a key role in our further considerations. In the following, let $x = (\xi_n)_{n=1}^\infty$ be a non-increasing bounded sequence of non-negative real numbers. For given $0<r< \infty$, $0 < s \le \infty$ and $\a \in \R$ we say that $x$ belongs to the generalized Lorentz sequence space $l_{r,s,\a}$ if 
$$ \left( (\log(n+1))^{-\a} \, n^{1/r-1/s} \, \xi_n \right)_{n=1}^\infty \in l_s, $$
i.e.
$$(\xi_n)_{n=1}^\infty \in l_{r,s,\a} \quad \text{if and only if} \quad \begin{cases} \sum\limits_{n=1}^\infty (\log(n+1))^{-\a s} \, n^{s/r-1} \, \xi_n^s < \infty , & \quad 0<s<\infty, \\ \sup\limits_{n \in \N} \, (\log(n+1))^{-\a} \, n^{1/r} \, \xi_n < \infty, & \quad s=\infty. \end{cases} $$
Note that we get the classical Lorentz sequence space $l_{r,s}$ for $\a=0$. A detailed explanation of classical Lorentz sequence spaces can be found, for instance, in \cite[Section 2.1]{P87}. If $A \subset M$ is a precompact subset of a metric space $M$, then the belonging of $(e_n(A))_n$ to a generalized Lorentz sequence space can be characterized by integral conditions of the \textit{metric entropy} $H(A,\e) := \log_2(N(A,\e))$ as follows: For $0<r<\infty$, $0 < s \le \infty$ and $\a \in \R$ we have that $(e_n(A))_n \in l_{r,s,\a}$ if and only if 
\begin{eqnarray} \label{rf1} 
\int_0^{e_1(A)} \left( \log_2(2+H(A,\e)) \right)^{-\a s} \, \left( H(A,\e) \right)^{s/r} \, \mathrm{d}\e^s < \infty 
\end{eqnarray}
for  $0 < s < \infty$ and
\begin{eqnarray} \label{rf2}
 \sup_{\e >0} \e \, \left( \log_2(2+H(A,\e)) \right)^{-\a} \, \left( H(A,\e) \right)^{1/r} < \infty 
\end{eqnarray}
for $ s = \infty$, respectively. This fact can be found in \cite[Lemma 2.2.]{CHP11}. Note that compared with \cite{CHP11} we use a slightly modified definition of generalized Lorentz sequence spaces. 

In general, the absolutely convex hull $\aco(A)$ of a subset $A \subset X$ of a Banach space $X$ is much larger than $A$ itself. Nevertheless, if $A$ is precompact, then it is well-known that also $\aco(A)$ is precompact. Hence, it seems natural to ask how the rate of decay of the entropy numbers of $A$ reflects the rate of decay of the entropy numbers of $\aco(A)$. In recent years this problem has been intensively studied in different settings (cf. e.g. \cite{BP90, C82, C85, C97, CE01, CE03, CHK88, CHP11, CKP99, CS90, CrSt02, D67, D73, D87, G01, G04, G12, GKS87, Kl12b, Kü05, Ky00, LL00, Pi81, St00, St04, Ta93}). In our setting, we consider a precompact subset $A$ of a Banach space of type $p$, $1 < p \le 2$. We assume that the dyadic entropy numbers of $A$ belong to a generalized Lorentz sequence space, which implies certain decay or summability properties. In Section \ref{chEn}, we explore the decay or summability properties of the dyadic entropy numbers of $\aco(A)$ expressed in terms of Lorentz spaces under the above-mentioned conditions. For further information about the relation of our results to the literature we refer to the introduction in Section \ref{chEn} and to \cite{CHP11}. Using (\ref{rf1}) and (\ref{rf2}), our results can easily be reformulated in terms of integral conditions of the metric entropy. 

The entropy of convex hulls and related topics have useful applications in analysis, approximation theory, geometry as well as probability (cf. e.g. \cite{CS90, ET96, LT91, Pi89}). In all these applications one is interested in sharp asymptotic estimates of entropy numbers. The study of convex hulls of precompact sets $A$ in Banach spaces is universal insofar as the entropy numbers $\e_n(S)$ of a compact operator $S:X \to Y$ between Banach spaces $X$ and $Y$ are always shared by the entropy numbers of a compact operator $T: l_1(A) \to Y$ defined on a $l_1(A)$ space in the sense that 
$$ \e_n(S) = \e_n(T). $$
This is one reason for studying the asymptotic behavior of entropy numbers of absolutely convex hulls. 

In Section \ref{chAp}, we will treat diverse applications of the entropy of convex hulls. We get new insight into the entropy of operators factoring through diagonal operators in $l_1$. Furthermore, we give optimal entropy estimates of (abstract) integral operators $T:X \to C(M)$ mapping from a Banach space $X$ into the space $C(M)$ of continuous functions on a compact metric space $(M,d)$. Our study is motivated by the universality of the Banach space $C[a,b]$ of continuous functions on the closed, bounded interval $[a,b]$ in the sense that any separable Banach space is isometrically isomorphic to a subspace of $C[a,b]$ (cf. \cite{B32}). As pointed out in \cite[Section 5.1]{CS90}, this implies universality of the class of $C[a,b]$-valued operators in the following sense: given a compact operator $T: X \to Y$ between arbitrary Banach spaces $X$ and $Y$, there is a compact operator $S: X \to C[a,b]$ such that the dyadic entropy numbers of $T$ are equivalent to those of $S$ in the sense that
$$ \frac{1}{2} \, e_n(S: X \to C[a,b]) \le e_n(T:X \to Y) \le 2 \, e_n(S:X \to C[a,b]). $$ 
In this sense, compact $C[a,b]$-valued operators represent a model for general compact operators. We shall see how the geometry of the Banach space $X$, the entropy numbers $\e_n(M)$ of the underlying compact metric space $M$ and the smoothness of the operator $T:X \to C(M)$ in terms of the \textit{modulus of continuity $\omega(T,\delta)$ of $T$},
$$ \omega(T,\delta) := \sup_{x \in B_X} \sup\left\{ | (Tx)(s) - (Tx)(t) | : s,t \in M, d(s,t) \le \delta \right\} $$
all affect the estimates of the dyadic entropy numbers of $T$. Finally, we apply these results to weakly singular integral operators mapping from an $L_p[0,1]$ space into the space $C[0,1]$ of continuous functions on the interval $[0,1]$ and into $L_q[0,1]$. In particular, we consider the classical Riemann-Liouville operator. During the last decades, a considerable amount of research has been devoted to the study of $s$-numbers, entropy numbers and eigenvalues of weakly singular operators and Riemann-Liouville operators in particular (cf. e.g. \cite{Bu07,CoKü88,CoKü90,DM97,Do93,Do95,FM86,Lif10,Lin04,Lin08,LL99,LS06,M01}).

For our purpose, we need some preliminaries. Unless otherwise stated, let $X$ denote a Banach space and $X'$ its dual Banach space. A Banach space $X$ is said to be \textit{of type $p$}, $1 \le p \le 2$, if there exists a constant $\rho > 0$ such that for all $k \in \N$ and all $x_1,\ldots,x_k \in X$ holds 	
$$ \int\limits_0^1 \left\Vert \sum_{i=1}^k r_i(t) \, x_i \right\Vert \, \mathrm{d}t \le \rho \left( \sum_{i=1}^k ||x_i||^p \right)^{1/p} . $$
Here $r_i$ denotes the $i$-th \textit{Rademacher function} given by
$$ r_i:[0,1] \to \R, \quad r_i(t) := \sgn(\sin(2^i\pi t)), $$ 
for $i=1,2,3,\ldots$. The \textit{type $p$ constant $\tau_p(X)$} is the smallest constant $\rho$ satisfying the above inequality. We say that $X$ has \textit{optimal type $p$}, if $X$ is of type $p$ but not of type $q$ for $q > p$. 
Furthermore, a Banach space is called \textit{B-convex}, if it is of some type $p > 1$. Recall that a Hilbert space has type $2$ and, due to the principle of local reflexivity, a Banach space and its bidual have the same type. Moreover, a Banach space is B-convex if and only if its dual is (cf. \cite{Pi73a,Pi73b} \cite[Corollary 13.7]{DJT95}). We are interested in the following \textit{s-numbers} associated with an operator $T \in \mathcal L(X,Y)$ (cf. \cite{P87}, \cite{CS90}): 

\begin{itemize}
	\item the \textit{$n$-th approximation number of $T$}, defined by
	$$ a_n(T) := \inf\left\{ ||T-A|| : A \in \mathcal L(X,Y) \text{ with } \rank A < n \right\}, $$
	\item the \textit{$n$-th Gelfand number of $T$}, defined by 
	$$ c_n(T) := \inf \left\{ ||T I_{E}^X || : E \text{ subspace of } X \text{ with } \codim(E) < n  \right\}, $$
	where $I_{E}^X$ is the natural embedding of $E$ into $X$,
	\item the \textit{$n$-th Kolmogorov number of $T$}, defined by
	$$ d_n(T) := \inf \left\{ ||Q_{F}^Y T|| : F \text{ subspace of } Y \text{ with } \dim(F) < n  \right\},   $$
	where $Q_{F}^Y: Y \to Y/F$ is the usual quotient map,
	\item the \textit{$n$-th symmetrized approximation number of $T$}, defined by 
	$$ t_n(T) := a_n(J_Y T Q_X), $$  
	where $Q_X$ is the canonical quotient map from $l_1(B_X)$ onto $X$ and $J_Y$ is the canonical embedding of $Y$ into $l_\infty(B_{Y'})$ (cf. \cite[p. 52, 60]{CS90}).
\end{itemize}
Recall that $t_n(T)=c_n(T Q_X)=d_n(J_Y T)$, $t_n(T) \le c_n(T), d_n(T) \le a_n(T)$ and $t_n(T)=t_n(T')$ holds true for all operators $T \in \mathcal L(X,Y)$, where $T'$ is the dual operator of $T$ (cf. \cite{P74,CS90}). The following inequality relates the entropy numbers of an operator to the above-mentioned $s$-numbers (cf. \cite[Theorem 1]{C81a}, \cite[Theorem 3.1.1]{CS90}, see also \cite[Theorem 1.3]{CKP99}). 

\begin{theo} \label{thC1}
For every $0<\a<\infty$ there exists a constant $c(\a) \ge 1$ such that for every operator $T \in \mathcal L(X,Y)$ between arbitrary Banach spaces $X$ and $Y$ and all $n \in \N$ we have 
$$ \sup_{1 \le k \le n} k^\a e_k(T) \le c(\a) \sup_{1 \le k \le n} k^\a t_k(T). $$ 
\end{theo}

Furthermore, we would like to point out that the entropy numbers $e_n(\aco(A))$ of the absolutely convex hull of a 
bounded subset $A$ of a Banach space $X$ can be expressed in terms of entropy 
numbers of operators. 
For this purpose, let $T_A : l_1(A) \to X$ be the operator defined by 
$$ T_A(a) := \sum_{t \in A} \a_t \, t ,\quad a=(\a_t)_{t \in A} \in l_1(A), $$
where $l_1(I)$ denotes the Banach space of all summable families $(\xi_t)_{t \in I}$ of real or complex numbers over the index set $I$, equipped with the norm
$$ || (\xi_t)_{t \in I} || = \sum_{t \in I} |\xi_t|. $$ 
Then 
\begin{align} \label{glaco}
\aco(A) \subset T_A(B_{l_1(A)}) \subset \overline{\aco(A)}
\end{align}
and, therefore, $e_n(\aco(A)) = e_n(T_A)$. 
Furthermore, the \textit{$n$-th Gelfand number of $\aco(A)$} is defined by
\begin{align} \label{glS4} 
c_n(\aco(A)) := c_n(T_A: l_1(A) \to X).  
\end{align}
For a geometrical interpretation of Gelfand numbers and a comparison with Gelfand widths we refer to the concluding remarks.

A subset $A$ of a vector space is called \textit{symmetric} if $A = -A$. The Banach space of all compact operators from $X$ into $Y$ is denoted by $\mathcal K(X,Y)$. For $1 \le p \le \infty$ we define the conjugate exponent $p'$ by the Hölder condition $1/p+1/p'=1$. By $\log(x)$ we denote the binary logarithm, i.e. the logarithm to the base two. Moreover, for $x \in \R$, the largest integer not greater than $x$ is denoted by $\left\lfloor x \right\rfloor := \max\{ z \in \Z \mid z \le x \}$. Note that $\left\lfloor x \right\rfloor \le x < \left\lfloor x \right\rfloor + 1 $ holds true for all $x \in \R$. In order to compare sequences, we introduce the following notations: For given sequences $(x_n)_n$, $(y_n)_n$ of positive numbers we write $x_n \preccurlyeq y_n$, if there exists a constant $c>0$ such that $x_n \le c \, y_n$ for all $n \in \N$. Furthermore, $x_n \sim y_n$ means that both $x_n \preccurlyeq y_n$ and $y_n \preccurlyeq x_n$.

\section{Entropy of absolutely convex hulls in Banach spaces of type $p$} \label{chEn}

Given a Banach space $X$ of type $p$, $1 < p \le 2$, and a precompact subset $A \subset X$ of $X$ it was shown in \cite[Proposition 6.2.]{CKP99} that
\begin{align} \label{glS1}
 (e_n(A))_n \in l_{r,\infty} \quad \text{implies} \quad (e_n(\aco(A)))_n \in l_{p',\infty,\a} 
\end{align}
holds true for all $r$ with $ 0 < r < p'$, where $\a=1/p'-1/r$. In the Hilbert space case, the corresponding result was proved in 
\cite[Proposition 5.5.]{CKP99}. Note that (\ref{glS1}) was extended by Steinwart who proved that, for all $r$ with $ 0< r < p'$ and $\gamma \in \R$,
$(e_n(A))_n \in l_{r,\infty,\gamma}$ implies
\begin{align} \label{glS5}
 e_n(\aco(A)) \preccurlyeq n^{-1/p'} \, (\log(n+1))^{1/p'-1/r} \, (\log\log(n+3))^\gamma
\end{align} 
(cf. \cite[Corollary 4]{St00}). We are interested in the setting of (\ref{glS1}) but for secondary index $s \not= \infty$, i.e. we ask for entropy estimates of $\aco(A)$ in terms of Lorentz sequence spaces under the condition that $(e_n(A))_n \in l_{r,s}$ for $ 0 < r < p'$ and $ 0 < s < \infty$.
In this context, Steinwart published the following result (cf. \cite[Theorem 1.3.]{St04}): Let $ 0 < s < \infty$ and define $r$ by $1/r=1/p'+1/s$, then 
\begin{align} \label{glS2}
(e_n(A))_n \in l_{r,s} \quad \text{implies} \quad (e_n(\aco(A)))_n \in l_{p',s}. 
\end{align}
However, since the choice of $s$ fixes $r$, this result does not have the desired generality. The following theorem closes this gap. 

\begin{theo} \label{TT02} Let $X$ be a Banach space of type $p$, $1 < p \le 2$, $0< r < p'$ and $ 0 < s < \infty$. Then for all precompact subsets 
$A \subset X$ we have that
$$ (e_n(A))_n \in l_{r,s} \quad \text{implies} \quad (e_n(\aco(A)))_n \in l_{p',s,\a} ,$$
where $\a = \frac{1}{s} + \frac{1}{p'} - \frac{1}{r} $. The result is optimal in the following sense: If $\beta < \a$ then there exists a precomact subset $A \subset l_p$ such that $(e_n(A))_n \in l_{r,s}$ and $(e_n(\aco(A) ))_n \notin l_{p',s,\beta}$.
\end{theo} 

\noindent Observe that Steinwart's result (\ref{glS2}) is contained in Theorem \ref{TT02}. Moreover, (\ref{glS1}) can be considered as the limit case $s=\infty$. The proof of Theorem \ref{TT02} is based on the following striking result of Steinwart (cf. \cite[Theorem 1.1]{St04}), which estimates the entropy numbers of $\aco(A)$ in terms of finitely many entropy numbers of $A$. 

\begin{theo} \label{T01} { \em \cite{St04} } Let $X$ be a Banach space of type $p$, $1 < p \le 2$, and $0<t< \infty$. Then there exists a constant $c(t)>0$ such that for
all integers $n \ge 2$, all integers $\a_1 < \a_2 < \ldots < \a_n$ and all bounded symmetric subsets $A \subset X$ we have 
\begin{eqnarray*}
e_{2m}(\aco(A)) & \le & c(t) \, m^{-1/t-1/p'} \sup_{i \le \min\{ m^{1+t/p'}; \, \a_1 \}} i^{1/t} \e_i(A) \\ 
& & + 23 \, \tau_p(X) \, 2^{-n/p'} \left( \sum_{k=1}^n \left( 2^{k/p'} \sum_{i=k}^n \e_{\a_i}(A) \right)^p \right)^{1/p},
\end{eqnarray*}
where $m := \left\lfloor 2^{n+2} \sum\limits_{k=2}^n 2^{-k} \log \left( \frac{2^{k+2}\a_k}{2^n} + 3 \right) \right\rfloor + 2$.
\end{theo}

Now let us consider the case of slow logarithmic decay of $(\e_n(A))_n$, i.e. $(e_n(A))_n \in l_{r,\infty}$ for $p' < r < \infty$. In \cite[Proposition 6.4.]{CKP99} it was shown that
\begin{align} \label{glS3}
  (e_n(A))_n \in l_{r,\infty} \quad \text{if and only if} \quad (e_n(\aco(A)))_n \in l_{r,\infty} .
\end{align}
The corresponding result for the Hilbert space case was proved in \cite[Proposition 5.5.]{CKP99}. With the help of the following inequality, Steinwart was able to refine the result (\ref{glS3}) by establishing a finite inequality as stated in the following theorem (cf. \cite[Theorem 4.2]{St99}, \cite[Theorem 4]{St00}). 

\begin{theo} \label{th03} 
{ \em \cite{St00} } Let $X$ be a Banach space of type $p$, $1 < p \le 2$. Then for all $r$ with $p' < r < \infty$ there exists a constant $c=c(p,r) > 0$, such that for
all $n \in \N$ and all precompact subsets $A \subset X$ we have
$$ \sup_{1 \le k \le n} k^{1/r} e_k(\aco(A)) \le c \, c_A \, \sup_{1 \le k \le n} k^{1/r} e_k(A), $$
where $c_A := \frac{\sup_{x \in A} ||x|| }{\e_1(A)} $. If $X$ is a Hilbert space and $p=2$, this is also true for the Gelfand numbers $c_k(\aco(A))$. 
\end{theo}

\noindent It turns out, that in B-convex Banach spaces the subsets $A$ and $\aco(A)$ surprisingly have the same entropy behavior, whenever $(e_n(A))_n$ or $(e_n(\aco(A)))_n$ decrease slow enough. More precisely (cf. \cite[Corollary 3]{St00}), if $p' < r < \infty$ and $(a_n)_n$ is a
positive sequence such that $(n^{1/r} a_n)_n$ is monotone increasing, then 
$$ e_n(A) \preccurlyeq a_n \quad \text{if and only if} \quad e_n(\aco(A)) \preccurlyeq a_n $$
and
$$ e_n(A) \sim a_n \quad \text{if and only if} \quad e_n(\aco(A)) \sim a_n .$$
The following theorem gives a positive answer to the open question whether (\ref{glS3}) holds true for secondary index $s \not= \infty$. 

\begin{theo} \label{TT03} Let $X$ be a Banach space of type $p$, $1 < p \le 2$, and let $p' < r < \infty$, $ 0 < s < \infty$ and $\a \in \R$. Then 
there exists a constant $c= c(p,r,s,\a) > 0$ such that for all $N \in \N$ and all precompact subsets $A \subset X$ it holds that
$$ \sum_{n=1}^N (\log(n+1))^{-\a} \, n^{s/r-1} (e_n(\aco(A)))^s \le c \, c_A^s \, \sum_{n=1}^N (\log(n+1))^{-\a} \, n^{s/r-1} (e_n(A))^s $$
and
$$ \sup_{1 \le n \le N} (\log(n+1))^{-\a} \, n^{1/r} e_n(\aco(A)) \le c \, c_A \, \sup_{1 \le n \le N} (\log(n+1))^{-\a} \, n^{1/r} e_n(A), $$
where $c_A := \frac{\sup_{x \in A} ||x|| }{\e_1(A)} $. In the context of Lorentz sequence spaces this implies that for $p' < r < \infty$, $0 < s \le \infty$ and $\a \in \R$ it holds that
$$ (e_n(A))_n \in l_{r,s,\a} \quad \text{if and only if} \quad (e_n(\aco(A)))_n \in l_{r,s,\a}. $$
\end{theo} 

What is left open is the so-called critical case of logarithmic decay of $(\e_n(A))_n$, i.e. $(e_n(A))_n \in l_{p',s}$ for $0 < s \le \infty$. For
$s=\infty$, Creutzig and Steinwart showed that 
\begin{align} \label{glS6}
 e_n(A) \le n^{-1/p'} \, (\log(n+1))^{-\gamma} \quad \text{implies} \quad e_n(\aco(A)) \preccurlyeq n^{-1/p'} \, (\log(n+1))^{-\gamma + 1} 
\end{align}
for $-\infty < \gamma < 1$. The result is asymptotically optimal (cf. \cite[Corollary 1.4.]{CrSt02}, see \cite{G01} for the Hilbert space case). The case $\gamma \ge 1$ is an open problem. 
In the case where $0 < s < \infty$, we conjecture that 
$$ (e_n(\aco(A)))_n \in l_{p',\infty,\b} \quad \text{with } \b = \max\left\{ 0 , 1-1/s \right\}. $$
However, we can not prove this and leave it as an open problem. Note that in the Hilbert space case, where $p=p'=2$, this problem has already been solved (cf. \cite{CE03}, \cite{CHP11}).

\section{Preliminary tools}

In this section we present basic tools for the proofs of the results. We start with the following inequalities of Hardy-type. 

\begin{lemma} \label{LH1}
Let $0<t<r<\infty$, $0<s<\infty$ and $\a \in \R$. If $\sigma_1 \ge \sigma_2 \ge \ldots \ge 0$ is a non-increasing sequence of non-negative real numbers then 
$$ \sum_{n=1}^N (\log(n+1))^\a \, n^{s/r-1} \, \left( \frac{1}{n} \sum_{k=1}^n \sigma_k^t \right)^{s/t} \le c \, \sum_{n=1}^N (\log(n+1))^\a \, n^{s/r-1} \, \sigma_n^s $$
for $N \in \N$, where $c = c(r,s,\a,t)>0$ is a constant depending on $r,s,\a$ and $t$.
\end{lemma}

A proof of this result can be found, for instance, in \cite[Lemma 2.3]{CHP11}. 

\begin{lemma} \label{LH2}
Let $0<t<r<\infty$ and $\a \in \R$. If $\sigma_1 \ge \sigma_2 \ge \ldots \ge 0$ is a non-increasing sequence of non-negative real numbers then 
$$ \sup_{1 \le n \le N} (\log(n+1))^\a \, n^{1/r} \left( \frac{1}{n} \sum_{k=1}^n \sigma_k^t \right)^{1/t} \le c \, \sup_{1 \le k \le N} (\log(k+1))^\a \, k^{1/r} \, \sigma_k $$
for $N \in \N$, where $c=c(t,r,\a)>0$ is a constant depending on $t,r$ and $\a$.
\end{lemma}

\proof Obviously we have 
\begin{eqnarray*}
\sum_{k=1}^n \sigma_k^t & = & \sum_{k=1}^n (\log(k+1))^{-\a t} \, k^{-t/r} \, (\log(k+1))^{\a t} \  k^{t/r} \, \sigma_k^t \\
& \le & \sup_{1 \le k \le n} (\log(k+1))^{\a t} \  k^{t/r} \, \sigma_k^t \, \sum_{k=1}^n (\log(k+1))^{-\a t} \, k^{-t/r} .
\end{eqnarray*}
Since $t/r<1$, we obtain
$$ \sum_{k=1}^n \sigma_k^t \le c(t,r,\a) \, (\log(n+1))^{-\a t} \, n^{-t/r+1} \sup_{1 \le k \le n} (\log(k+1))^{\a t} \  k^{t/r} \, \sigma_k^t,  $$
which yields
$$ (\log(n+1))^{\a t} \, n^{t/r} \, \frac{1}{n} \sum_{k=1}^n \sigma_k^t \le c(t,r,\a) \, \sup_{1 \le k \le n} (\log(k+1))^{\a t} \  k^{t/r} \, \sigma_k^t. $$
Hence, for all $n \in \N$, we get
\begin{eqnarray*}
(\log(n+1))^{\a} \, n^{1/r} \, \left( \frac{1}{n} \sum_{k=1}^n \sigma_k^t \right)^{1/t} & \le & c(t,r,\a) \left( \sup_{1 \le k \le n} (\log(k+1))^{\a t} \  k^{t/r} \, \sigma_k^t \right)^{1/t} \\
& = & c(t,r,\a) \sup_{1 \le k \le n} (\log(k+1))^\a \, k^{1/r} \, \sigma_k
\end{eqnarray*}
and taking the supremum with respect to $1 \le n \le N$ finishes the proof. \hfill $\blacksquare$

In order to prove the optimality of the result stated in Theorem \ref{TT02} we will use the following lemma (cf. \cite{CKP99}, \cite[Lemma 2.4]{CHP11}).

\begin{lemma} \label{L02}
Let $\sigma_1 \ge \sigma_2 \ge \ldots \ge 0$ be a non-increasing sequence of non-negative real numbers. For $1 < p \le 2$, let 
$$ A = \{ \sigma_n u_n \mid n \in \N \} \subset l_p, $$
where $\{u_1,u_2,\ldots\}$ denotes the canonical unit vector basis of the sequence space $l_p$. Then for all $n \in \N$ we have
$$ \e_n(A) \le \sigma_n \quad \text{and} \quad e_n(\aco(A)) \ge c \, \max\{ n^{-1/p'} (\log(n+1))^{1/p'} \sigma_{n^2} , \sigma_{2^n} \}, $$
where $c=c(p)>0$ is a constant only depending on $p$.
\end{lemma}

\proof Since $|| \sigma_k u_k \mid l_p || = \sigma_k \le \sigma_n $ for all $k \ge n$ it is obvious that $\e_n(A) \le \sigma_n$. To estimate the entropy numbers of the absolutely convex hull of $A$, we consider the sections 
$$ \Delta_{n,m} := \aco \{ \sigma_k u_k \mid n \le k \le m \} , \quad m,n \in \N, m>n.  $$ 
Due to the monotonicity we have
$$ e_n(\aco(A)) \ge e_n(\Delta_{n,m}) \ge \sigma_m \, e_n(id:l_1^{m-n} \to l_p^{m-n}) $$
and by a result of Schütt \cite{Sch84} and Garnaev/Gluskin \cite{GG84} (see also \cite{CP88} for a generalization) it holds
$$ e_n(id:l_1^{m-n} \to l_p^{m-n}) \ge c \, \left( \frac{\log(m/n)}{n}  \right)^{1/p'}, $$
where $c>0$ is an absolute constant. By putting $m=n^2$ and $m=2^n$, respectively, the assertion follows. \hfill $\blacksquare$

Next, we recall a striking result of Bourgain et al., which relates the entropy numbers of a compact operator to the entropy numbers of its dual operator (cf. \cite{BPST89}, see \cite{To87} for the Hilbert space case).

\begin{theo} \label{endu}
Let $X$ and $Y$ be Banach spaces such that one of them is B-convex. Then for every $0< \a <\infty$ there exists a constant $c=c(\a) \ge 1$ such that for all compact operators $T: X \to Y$ and all $n \in \N$ we have
$$ c^{-1} \, \sup_{1 \le k \le n} k^{\a}e_k(T) \le \sup_{1 \le k \le n} k^{\a}e_k(T') \le c \, \sup_{1 \le k \le n} k^{\a}e_k(T). $$
\end{theo}

There is an improved version of this duality result of Artstein et al. in \cite{AMST04} which we could have used instead.
But Theorem \ref{endu} is sufficient for our purpose.

\section{Proof of the results}

\noindent \textit{Proof of Theorem \ref{TT02}.} The proof uses techniques and ideas from \cite{St04}. It is enough to prove the result for symmetric subsets $A$. We choose a constant $a$ with $\frac{1}{2} < a < 1$ and define $\a_k := \left\lfloor 2^{n2^{a(k-1)}} \right\rfloor$ for $k=1,2,\ldots,n$, then $ 2^n = \a_1 < \a_2 < \ldots < \a_n $. Furthermore, with 
$$ m = 2 +  \left\lfloor 2^{n+2} \sum\limits_{k=2}^n 2^{-k} \log \left( \frac{2^{k+2} \left\lfloor 2^{n2^{a(k-1)}} \right\rfloor }{2^n} + 3 \right) \right\rfloor $$ 
it holds 
\begin{eqnarray*}
c_1 \, n 2^n  \le m \le c_2 \, n 2^n, 
\end{eqnarray*}
where $c_1,c_2 > 0$ are constants depending on $a$. Since $2^{\left\lfloor n 2^{a(i-1)} \right\rfloor} \le \left\lfloor 2^{n2^{a(i-1)}} \right\rfloor $ for $i \in \N$, we have that $\e_{\a_i}(A) \le e_{\left\lfloor n 2^{a(i-1)} \right\rfloor+1}(A) $. Consequently, we obtain 
\begin{eqnarray*}
\left( \sum_{k=1}^n \left( 2^{k/p'} \sum_{i=k}^n \e_{\a_i}(A) \right)^p \right)^{1/p} & \le & \sum_{k=1}^n 2^{k/p'} \sum_{i=k}^n e_{\left\lfloor n 2^{a(i-1)} \right\rfloor+1}(A) \\ 
& = &  \sum_{i=1}^n e_{\left\lfloor n 2^{a(i-1)} \right\rfloor+1}(A) \sum_{k=1}^i 2^{k/p'} \\ 
& \le & c_3 \sum_{i=1}^n 2^{i/p'} e_{\left\lfloor n 2^{a(i-1)} \right\rfloor+1}(A),
\end{eqnarray*}
where $c_3 > 0$ only depends on $p$. Hence, Theorem \ref{T01} yields
$$ \begin{aligned}
e_{2c_2n2^n}(\aco(A)) \le & \; e_{2m}(\aco(A)) \le c_4 \left( n2^n \right)^{-1/t-1/p'} \sup_{1 \le i \le 2^n} i^{1/t} \e_i(A) \\ & + c_5 \, 2^{-n/p'} \sum_{i=1}^n 2^{i/p'} e_{\left\lfloor n 2^{a(i-1)} \right\rfloor+1}(A)
\end{aligned} $$
for all $n \ge 2$ and all $t>0$, where $c_4,c_5 > 0$ do not depend on $n$ and $A$. By standard arguments, we conclude that
$$ e_{2c_22^n}(\aco(A)) \le c_6 \, 2^{-n/t-n/p'} \sup\limits_{1 \le i \le 2^n} i^{1/t} \e_i(A) + c_7 \, n^{1/p'} 2^{-n/p'} \sum_{i=1}^n 2^{i/p'} e_{\left\lfloor n 2^{a(i-1)-1} \right\rfloor+1}(A) $$
for all $n \ge 2$ and all $t > 0$, where $c_6,c_7 > 0$ are constants independent of $n$ and $A$. Using a dyadic characterization of Lorentz sequence
spaces (cf. \cite[2.1.10.*]{P87}), we see that the assertion is equivalent to 
$$ \left( (n+1)^{-\a} \, 2^{n/p'} \, e_{2^n}(\aco(A)) \right)_n \in l_s , \quad \a = \frac{1}{s} + \frac{1}{p'} - \frac{1}{r}, $$
hence, it suffices to show that 
\begin{enumerate}
\item[(1)] $ \left( n^{-\a} \, 2^{-n/t} \sup\limits_{1 \le i \le 2^n} i^{1/t} \e_i(A) \right)_n \in l_s $
\item[(2)] $ \left( n^{-\a+1/p'} \sum\limits_{i=1}^n 2^{i/p'} e_{\left\lfloor n 2^{a(i-1)-1} \right\rfloor+1}(A) \right)_n \in l_s $
\end{enumerate}
for a suitable $t > 0$. First let us deal with (1). Since $(e_n(A))_n \in l_{r,s} \subset l_{r,\infty}$, there exists a constant $c_8 >0$ such that
$\e_n(A) \le c_8 \, (\log(n+1))^{-1/r}$ for all $n\in \N$. We fix $t$ with $0 < t \le r$ and obtain
$$ \sup\limits_{1 \le i \le 2^n} i^{1/t} \e_i(A) \le c_8 \sup\limits_{1 \le i \le 2^n} \frac{i^{1/t}}{(\log(i+1))^{1/r}} \le c_9 \, n^{-1/r} \, 2^{n/t}. $$  
Hence, 
$$ n^{-\a} \, 2^{-n/t} \sup\limits_{1 \le i \le 2^n} i^{1/t} \e_i(A) \le c_9 \, n^{-\a-1/r} = c_9 \, n^{-1/p'-1/s}, $$
which yields
$$ \sum_n \left( n^{-\a} \, 2^{-n/t} \sup\limits_{1\le i \le 2^n} i^{1/t} \e_i(A) \right)^s \le c_9^s \sum_n n^{-s/p'-1} < \infty. $$ 
The proof of (2) is more technical. Let $\e >0$ be arbitrary. First of all we show that 
$$ \left( \sum\limits_{i=1}^n 2^{i/p'} e_{\left\lfloor n 2^{a(i-1)-1} \right\rfloor+1}(A) \right)^s \le c_{10} \sum\limits_{i=1}^n 2^{ib} \left(e_{\left\lfloor n 2^{a(i-1)-1} \right\rfloor+1}(A) \right)^s, $$
where $b := (1+\e)\frac{s}{p'}$ and $c_{10} > 0$ only depends on $s$ and $p$. In the case $0<s \le 1$ this is obvious since 
$$ \left( \sum\limits_{i=1}^n 2^{i/p'} e_{\left\lfloor n 2^{a(i-1)-1} \right\rfloor+1}(A) \right)^s \le \sum_{i=1}^n 2^{is/p'} \left(e_{\left\lfloor n 2^{a(i-1)-1} \right\rfloor+1}(A) \right)^s .$$
If $1<s<\infty$ we use Hölder's inequality to see that
\begin{eqnarray*}
& &\left( \sum_{i=1}^n 2^{i/p'} e_{\left\lfloor n 2^{a(i-1)-1} \right\rfloor+1}(A) \right)^s \\
& = &  \left( \sum_{i=1}^n \left( 2^{is/p'-ib} \, 2^{ib} \left( e_{\left\lfloor n 2^{a(i-1)-1} \right\rfloor+1}(A) \right)^s \right)^{1/s} \right)^s \\
& \le & \left( \sum_{i=1}^n 2^{(is/p'-ib)/(s-1)} \right)^{s-1} \sum_{i=1}^n 2^{ib} \left( e_{\left\lfloor n 2^{a(i-1)-1} \right\rfloor+1}(A) \right)^s .
\end{eqnarray*} 
Because of $(s/p' - b)/(s-1) < 0$ we obtain the desired result. Consequently, we get 
\begin{eqnarray} \label{glt1}
& & \sum_{n=1}^N \left( n^{-\a+1/p'} \sum\limits_{i=1}^n 2^{i/p'} e_{\left\lfloor n 2^{a(i-1)-1} \right\rfloor+1}(A) \right)^s \nonumber \\
& \le & c_{10} \sum_{n=1}^N n^{s/r-1} \sum\limits_{i=1}^n 2^{ib} \left(e_{\left\lfloor n 2^{a(i-1)-1} \right\rfloor+1}(A) \right)^s \\
& = & c_{10} \sum_{i=1}^N 2^{ib} \sum\limits_{n=i}^N n^{s/r-1} \left(e_{\left\lfloor n 2^{a(i-1)-1} \right\rfloor+1}(A) \right)^s . \nonumber
\end{eqnarray}
In a last step we check that 
\begin{align} \label{glt2}
 \sum_{i=1}^N 2^{ib} \sum\limits_{n=i}^N n^{s/r-1} \left(e_{\left\lfloor n 2^{a(i-1)-1} \right\rfloor+1}(A) \right)^s \le c_{11} \sum_{i=1}^N 2^{ib-ias/r} \sum_{n=1}^\infty n^{s/r-1} (e_n(A))^s 
\end{align}
with a constant $c_{11}>0$ independent of $A$ and $N$. For the sake of simplicity we introduce the notation $a_i := a(i-1)-1$. First observe that there exists a constant $C_1(r,s)>0$  such that for all $n,i \in \N$ with $n \ge i$ it holds $\left( n2^{a_i} \right)^{s/r-1} \le C_1(r,s) \, \left( \left\lfloor n2^{a_i} \right\rfloor +1 \right)^{s/r-1} $. Hence, the estimate
\begin{eqnarray*}
\sum_{n=i}^N n^{s/r-1} \left( e_{ \left\lfloor n2^{a_i} \right\rfloor + 1 }(A) \right)^s & = & 2^{-a_i(s/r-1)} \sum_{n=i}^N \left( n2^{a_i} \right)^{s/r-1} \left( e_{ \left\lfloor n2^{a_i} \right\rfloor + 1 }(A) \right)^s \\
& \le & C_2(r,s) \, 2^{-ia(s/r-1)} \sum_{n=i}^N \left( \left\lfloor n2^{a_i} \right\rfloor + 1 \right)^{s/r-1} \left( e_{ \left\lfloor n2^{a_i} \right\rfloor + 1 }(A) \right)^s 
\end{eqnarray*}
hold. Next, due to the monotonicity, we have
$$ \left( e_{\left\lfloor n2^{a_i} \right\rfloor+1}(A) \right)^s \le \frac{ \left( e_{\left\lfloor (n-1)2^{a_i} \right\rfloor+1}(A) \right)^s + \left( e_{\left\lfloor (n-1)2^{a_i} \right\rfloor+2}(A) \right)^s + \ldots + \left( e_{\left\lfloor n2^{a_i} \right\rfloor+1}(A) \right)^s  }{ \left\lfloor n2^{a_i} \right\rfloor-\left\lfloor (n-1)2^{a_i} \right\rfloor+1 } $$
and since
$$ \left\lfloor n2^{a_i} \right\rfloor-\left\lfloor (n-1)2^{a_i} \right\rfloor+1 \ge n2^{a_i}-(n-1)2^{a_i} = 2^{a_i} $$
we conclude
\begin{eqnarray*}
\left( e_{\left\lfloor n2^{a_i} \right\rfloor+1}(A) \right)^s & \le & 2^{-a_i} \left(  \left( e_{\left\lfloor n2^{a_i} \right\rfloor+1}(A) \right)^s + \sum_{k=\left\lfloor (n-1)2^{a_i} \right\rfloor + 1 }^{\left\lfloor n2^{a_i} \right\rfloor } \left(e_k(A)\right)^s \right) \\
& \le & 2^{-a_i+1} \sum_{k=\left\lfloor (n-1)2^{a_i} \right\rfloor + 1 }^{\left\lfloor n2^{a_i} \right\rfloor} \left(e_k(A)\right)^s \\
& \le & C_3 \, 2^{-ia} \sum_{k=\left\lfloor (n-1)2^{a_i} \right\rfloor + 1 }^{\left\lfloor n2^{a_i} \right\rfloor} \left(e_k(A)\right)^s,  
\end{eqnarray*}
where $C_3 > 0$ is an absolute constant. Consequently, it holds
\begin{eqnarray*}
& & \sum_{i=1}^N 2^{ib} \sum\limits_{n=i}^N n^{s/r-1} \left(e_{\left\lfloor n 2^{a(i-1)-1} \right\rfloor+1}(A) \right)^s \\
& \le & C_2(r,s) \sum_{i=1}^N 2^{ib-ia(s/r-1)} \sum_{n=i}^N \left( \left\lfloor n2^{a_i} \right\rfloor + 1 \right)^{s/r-1} \left( e_{ \left\lfloor n2^{a_i} \right\rfloor + 1 }(A) \right)^s \\
& \le & C_4(r,s) \sum_{i=1}^N 2^{ib-ia(s/r-1)-ia} \sum_{n=i}^N \sum_{k=\left\lfloor (n-1)2^{a_i} \right\rfloor + 1 }^{\left\lfloor n2^{a_i} \right\rfloor} \left( \left\lfloor n2^{a_i} \right\rfloor + 1 \right)^{s/r-1} \left(e_k(A)\right)^s.
\end{eqnarray*}
Now we claim that there exists a constant $C_5(r,s)>0$ such that for all $n, i \in \N$ with $n \ge i$ it holds $ \left( \left\lfloor n2^{a_i} \right\rfloor + 1 \right)^{s/r-1} \le C_5(r,s) \, k^{s/r-1}$. If $s/r-1<0$ then this is obvious since $k \le \left\lfloor n2^{a_i} \right\rfloor + 1$. To handle the case $s/r-1\ge 0$ we observe that 
$$ \frac{\left\lfloor n2^{a_i} \right\rfloor+1}{\left\lfloor (n-1)2^{a_i} \right\rfloor+1} \le \frac{n2^{a_i}+1}{(n-1)2^{a_i}} = \frac{n}{n-1} + \frac{1}{(n-1)2^{a_i}} \le 4 $$
for all $n \ge 2$ and all $i \in \N$. Therefore, we obtain 
\begin{equation} \label{gl01} 
4k \ge 4 \left( \left\lfloor (n-1)2^{a_i} \right\rfloor + 1 \right) \ge \left( \left\lfloor n2^{a_i} \right\rfloor + 1 \right)
\end{equation}
and hence the assertion for all $n \ge i \ge 2$. Thus, left open is the case $i=1$. Since $2^{a_1} = \frac{1}{2}$ we have to show that 
$$ \left( \left\lfloor n/2 \right\rfloor + 1 \right)^{s/r-1} \le C_6(r,s) \, k^{s/r-1}, \quad s/r-1 \ge 0, $$
for all $n \in \N$, where $ \left\lfloor (n-1)/2 \right\rfloor + 1 \le k \le \left\lfloor n/2 \right\rfloor $. However, this is clear because it holds $k \ge \left\lfloor (n-1)/2 \right\rfloor + 1 \ge \frac{1}{2} \left( \left\lfloor n/2 \right\rfloor + 1 \right) $ for all $n \in \N$. Hence we conclude
\begin{eqnarray*}
& & \sum_{i=1}^N 2^{ib} \sum\limits_{n=i}^N n^{s/r-1} \left(e_{\left\lfloor n 2^{a(i-1)-1} \right\rfloor+1}(A) \right)^s \\
& \le & C_7(r,s) \sum_{i=1}^N 2^{ib-ias/r} \sum_{n=i}^N \sum_{k=\left\lfloor (n-1)2^{a_i} \right\rfloor + 1 }^{\left\lfloor n2^{a_i} \right\rfloor} k^{s/r-1} \left(e_k(A)\right)^s \\
& \le & C_7(r,s) \sum_{i=1}^N 2^{ib-ias/r} \sum_{m=1}^\infty m^{s/r-1} \left( e_m(A) \right)^s ,
\end{eqnarray*}
which is the desired estimate.  

Therefore, combining (\ref{glt1}) and (\ref{glt2}) we find that
$$ \sum_{n=1}^N \left( n^{-\a+1/p'} \sum\limits_{i=1}^n 2^{i/p'} e_{\left\lfloor n 2^{a(i-1)-1} \right\rfloor+1}(A) \right)^s \le c_{12} \sum_{i=1}^N 2^{ib-ias/r} \sum_{n=1}^\infty n^{s/r-1} (e_n(A))^s $$
for a constant $c_{12}>0$ independent of $A$ and $N$. Remember that $\frac{1}{2} < a < 1$ and $b = (1+\e)\frac{s}{p'}$ with $\e >0$ arbitrary. Now fix $\e >0$ such that $(1+\e)\frac{r}{p'} <1$. Then we can choose $a$ as a constant satisfying both $\frac{1}{2} < a < 1$ and $a > (1+\e)\frac{r}{p'}$. Consequently we have $b-as/r <0$, which yields the statement.  

Finally we will prove the optimality of the result. Let $\beta < \a = \frac{1}{s} + \frac{1}{p'} - \frac{1}{r} $. Choose $\gamma$ with $\gamma s > 1$ and consider 
$$ A = \{ \sigma_n u_n \mid n \in \N \} \subset l_p, $$
where $\sigma_n = (\log(n+1))^{-1/r} (\log \log(n+3))^{-\gamma}$. Taking Lemma \ref{L02} into account, we obtain 
$$ e_n(A) \le \sigma_{2^{n-1}} \le C_8(r,s) \, n^{-1/r} \, (\log(n+1))^{-\gamma} $$
and hence it holds
$$ \sum_n n^{s/r-1} (e_n(A))^s \le C_8^s(r,s) \sum_n n^{-1} \, (\log(n+1))^{-\gamma s}. $$
Since the latter series is convergent, we see that $(e_n(A))_n \in l_{r,s}$. Moreover, Lemma \ref{L02} yields 
$$ e_n(\aco(A)) \ge C_9(p,r,s) \, n^{-1/p'} \, (\log(n+1))^{1/p'-1/r} \, (\log \log(n+3))^{-\gamma} $$
and consequently we have
\begin{eqnarray*}
& & \sum_n (\log(n+1))^{-\beta s} \,  n^{s/p'-1}  \, (e_n(\aco(A)))^s \\ 
& \ge & C_9^s(p,r,s) \sum_n n^{-1} \, (\log(n+1))^{s(1/p'-1/r-\b)} \, (\log \log(n+3))^{-\gamma s} \\ 
& = & C_9^s(p,r,s) \sum_n n^{-1} \, (\log(n+1))^{\a s - \b s - 1} \, (\log \log(n+3))^{-\gamma s} = \infty, 
\end{eqnarray*}
which means $(e_n(\aco(A)))_n \notin l_{p',s,\beta}$. \hfill $\blacksquare$ 

\noindent \textit{Proof of Theorem \ref{TT03}.} The proof is based on Theorem \ref{th03} in combination with the Hardy-type inequalities from Lemma \ref{LH1} and \ref{LH2}. Choose $t$ with $p' < t < r$. Then according to Theorem \ref{th03} it holds 
\begin{eqnarray*}
n^{1/t} e_n(\aco(A)) \le \sup_{1 \le k \le n} k^{1/t} e_k(\aco(A)) \le c_1 \, c_A \, \sup_{1 \le k \le n} k^{1/t} e_k(A), 
\end{eqnarray*}
where $c_1 > 0$ only depends on $p$ and $t$. In addition, due to the monotonicity of the entropy numbers we have
$$ \sup_{1 \le k \le n} k^{1/t} e_k(A) \le \left( \sum_{k=1}^n (e_k(A))^t \right)^{1/t} $$ 
which yields
$$ n^{1/t} e_n(\aco(A)) \le c_1 \, c_A \, \left( \sum_{k=1}^n (e_k(A))^t \right)^{1/t} $$
and hence 
$$ e_n(\aco(A)) \le c_1 \, c_A \, \left( \frac{1}{n} \sum_{k=1}^n (e_k(A))^t \right)^{1/t}. $$
Consequently, in the case $0<s<\infty$ we obtain
$$ \begin{aligned} 
\sum_{n=1}^N (\log(n+1))^{-\a} & \, n^{s/r-1} (e_n(\aco(A)))^s \\
& \le c_1^s \, c_A^s \sum_{n=1}^N (\log(n+1))^{-\a} \, n^{s/r-1} \left( \frac{1}{n} \sum_{k=1}^n (e_k(A))^t  \right)^{s/t} 
\end{aligned} $$
and since $t<r$, we can use Lemma \ref{LH1} to get 
$$ \sum_{n=1}^N (\log(n+1))^{-\a} \, n^{s/r-1} (e_n(\aco(A)))^s \le c \, c_A^s \sum_{n=1}^N (\log(n+1))^{-\a} \, n^{s/r-1} (e_n(A))^s,  $$  
where $c>0$ is a constant independent of $A$ and $N$. The second inequality can be treated analogously by using Lemma \ref{LH2}. \hfill $\blacksquare$

\section{Applications of entropy and Gelfand numbers of convex hulls} \label{chAp}

We start with studying the entropy numbers of operators factoring through diagonal operators in $l_1$ and with values in a Banach space of type $p$. After that, we deal with $C(M)$-valued operators $T_K$ generated by abstract kernels $K$. Using a general approach, we show how entropy and Kolmogorov numbers of such abstract kernel integration operators are connected to the entropy of the image of the abstract kernel. We consider the important case of a weakly singular integral operator $T_K : L_p[0,1] \to C[0,1]$ generated by a convolution kernel. Here, the Hilbert space case $p=2$ is of particular interest. In addition, we also investigate the case when $T_K: L_2[0,1] \to L_q[0,1]$. Finally, we deal with entropy estimates of the classical Riemann-Liouville operator in different settings. In all these applications we need sharp estimates of entropy and Gelfand numbers of absolutely convex hulls. 

\subsection{Operators factoring through diagonal operators in $l_1$}

First, we want to demonstrate how the results of Theorem \ref{TT02} and \ref{TT03} can be used to give entropy estimates of operators $T$ admitting a factorization $T=SD$ where $D: l_1 \to l_1$ is a diagonal operator and $S: l_1 \to X$ is an arbitrary operator from $l_1$ into a Banach space $X$ of type $p$, $1 < p \le 2$. This complements a result given in \cite[Theorem 1]{C82} by adding the case of logarithmic decay.

\begin{theo} \label{THdi} Let $X$ be a Banach space of type $p$, $1 < p \le 2$, and $S \in \mathcal L(l_1,X)$. Furthermore, let $D_\sigma: l_1 \to l_1$ be the diagonal operator induced by the non-increasing and non-negative sequence $\sigma=(\sigma_n)_n$. Then for the composition operator $T=SD$ the following statements hold: 

\begin{enumerate}
\item[(i)] If $0 < r < p'$ and $0 < s \le \infty$, then 
$$ (\sigma_{2^{n-1}})_n \in l_{r,s} \quad \text{implies} \quad (e_n(T))_n \in l_{p',s,\a} \quad \text{with } \a = \frac{1}{s}+\frac{1}{p'}-\frac{1}{r}. $$  
\item[(ii)] If $p' < r < \infty$, $0 < s \le \infty$ and $\b \in \R$, then 
$$ (\sigma_{2^{n-1}})_n \in l_{r,s,\b} \quad \text{implies} \quad (e_n(T))_n \in l_{r,s,\b}. $$  
\end{enumerate}  
\end{theo}

\proof Consider the set 
$$ A := \{ T u_n \mid n \in \N \} \subset X, $$
where $\{u_1,u_2,\ldots\}$ denotes the canonical unit vector basis of $l_1$. Then, obviously, for all $n \in \N$ we have $\e_n(A) \le ||S|| \, \sigma_n$. Consequently, 
$$ (\sigma_{2^{n-1}})_n \in l_{r,s,\b}  \quad \text{implies} \quad (e_n(A))_n \in l_{r,s,\b}.  $$    
Furthermore, it holds 
$$ \aco(A) \subset T(B_{l_1}) \subset \overline{\aco(A)} $$
and hence we can conclude that 
$$ e_n(T) = e_n(\aco(A)). $$
Therefore, by applying Theorem \ref{TT02} and \ref{TT03}, we get the desired results. \hfill $\blacksquare$

Note that the results of Theorem \ref{THdi} are optimal. To see this, choose $X=l_p$ with $1 < p \le 2$ and $S = id: l_1 \to l_p$. Then $ e_n(T:l_1 \to l_p) = e_n(\aco(A)) $, where
$$ A = \{ \sigma_n u_n \mid n \in \N \} \subset l_p. $$
The assertion follows from Lemma \ref{L02} using the same argumentation as in the proof of the optimality in Theorem \ref{TT02}.

\subsection{Operators with values in $C(M)$}

As already mentioned in the introduction, the entropy behavior of a compact operator is reflected by that of a $C(M)$-valued operator on a compact metric space $(M,d)$. Thus, for our purpose, $C(M)$-valued operators are universal. By the Arzelà-Ascoli theorem, we know that an operator $T: X \to C(M)$ from a Banach space $X$ into the space $C(M)$ of all continuous scalar-valued functions on a compact metric space $M$ is compact if and only if the limit relation
\begin{align} \label{gl02}
	\lim_{\delta \to 0+} \omega(T,\delta) = 0
\end{align} 
is fulfilled (cf. \cite[Proposition 5.5.1]{CS90}). In analogy to Hölder-continuous functions, a compact operator $T:X \to C(M)$ is called \textit{Hölder-continuous of type $\a$, $ 0 < \a \le 1$}, if  
$$ |T|_\a := \sup_{\delta >0} \frac{\omega(T,\delta)}{\delta^\a} < \infty. $$
Note that a Hölder-continuous operator $T$ of type $\a$ actually maps $X$ into the space $C^\a(M)$ of Hölder-continuous scalar-valued functions of type $\a$ on $M$. In the special case $\a=1$, a Hölder-continuous operator of type $1$ is said to be \textit{Lipschitz-continuous}. The vector space $\mathcal Lip_\a(X,C(M))$ of all operators from $X$ into $C(M)$ which are Hölder-continuous of type $\a$ becomes a Banach space under the norm 
$$ \Lip_\a(T) := \max \left\{ ||T:X \to C(M)||,|T|_\a  \right\}. $$
For the space of all Lipschitz-continous operators $T:X \to C(M)$ we simply write 
$$ [\mathcal Lip(X,C(M)),\Lip] := [\mathcal Lip_1(X,C(M)),\Lip_1]. $$
By changing the metric $d$ on $M$ to $d^\a$, $0< \a \le 1$, we reduce a Hölder-continuous operator of type $\a$ to a Lipschitz-continuous operator, i.e.
$$ \mathcal Lip_\a(X,C( (M,d) ) ) = \mathcal Lip(X,C( (M,d^\a) ) ). $$ 

Now we represent compact and Hölder-continuous $C(M)$-valued operators by abstract kernels. To this end, let us introduce the vector space $C(M,Z)$ of all continuous $Z$-valued functions on a compact metric space $(M,d)$, where $Z$ is an arbitrary Banach space. It is clear that $C(M,Z)$ is a Banach space with respect to the supremum norm
$$ ||K||_\infty := \sup_{s \in M} ||K(s)||_Z. $$  
Just as for scalar-valued functions on $M$ we define a modulus of continuity by 
$$ \omega_Z(K,\delta) := \sup \left\{ ||K(s)-K(t)||_Z : s,t \in M, d(s,t) \le \delta \right\} $$
for $0 \le \delta < \infty$. Since $ \omega_Z(K,\delta) \le 2 \, ||K||_\infty $ holds true for all $\delta \ge 0$, this is well defined for arbitrary bounded $Z$-valued functions $K$ on $M$. From the fact that a continuous function on a compact set is uniformly continuous, we see that a bounded $Z$-valued function is continuous if and only if 
\begin{align} \label{gl03}
   \lim_{\delta \to 0+} \omega_Z(K,\delta) = 0. 
\end{align}
A stronger condition than the limit relation in (\ref{gl03}) is that of \textit{Hölder-continuity}. The continuous function $K \in C(M,Z)$ is said to be \textit{Hölder-continuous of type $\a$, $0 < \a \le 1$}, if 
$$ |K|_{Z,\a} := \sup_{\delta > 0} \frac{\omega_Z(K,\delta)}{\delta^\a} < \infty .$$
As a direct consequence of this definition, we see that $\omega_Z(K,\delta) \le |K|_{Z,\a} \, \delta^\a$ holds for all $\delta \ge 0$. Consequently, for all $s,t \in M$ we have
$$ ||K(s)-K(t)||_Z \le |K|_{Z,\a} \, (d(s,t))^\a $$
and this implies
$$ \e_n(\image(K)) \le |K|_{Z,\a} \, (\e_n(M))^\a $$
for all $n \in \N$. The vector space $C^\a(M,Z)$ of all Hölder-continuous $Z$-valued functions of type $\a$ on $M$ turns out to be a Banach space with respect to the norm 
$$ ||K||_{Z,\a} := \max\left\{||K||_\infty, |K|_{Z,\a} \right\}. $$

If $Z$ is a dual space, i.e. $Z=X'$ for some Banach space $X$, then an element $K \in C(M,X')$ gives rise to an operator $T_K: X \to C(M) $ according to the rule 
\begin{align} \label{glTK} 
(T_K x)(s) := \left\langle x, K(s) \right\rangle, \quad x \in X, s \in M,
\end{align}
where $\left\langle \cdot, \cdot \right\rangle$ denotes the duality pairing between $X$ and $X'$. The function $K$ is called the \textit{abstract kernel} of the operator $T_K$. Obviously, $T_K$ is a linear operator and an easy computation shows that it is also bounded: 
$$ ||T_K : X \to C(M)|| = \sup_{x \in B_X} \sup_{s \in M} |\left\langle x, K(s) \right\rangle| = \sup_{s \in M} ||K(s)||_{X'} = ||K||_\infty < \infty. $$ 
Moreover, we have that 
\begin{align}
\omega(T_K,\delta) = \omega_{X'}(K,\delta), 
\end{align}
hence, by combining (\ref{gl02}) and (\ref{gl03}), we observe that the operator $T_K: X \to C(M)$ is compact. It does not surprise that stronger conditions on the kernel $K$ are reflected in the properties of the generated operator $T_K$. Indeed, if the kernel $K$ is even Hölder-continuous of type $\a$, then $T_K:X \to C(M)$ is a Hölder-continuous operator of type $\a$.  

On the other hand, any compact or even Hölder-continuous operator $T$ from $X$ into $C(M)$ can be generated by an appropriate kernel in the sense of (\ref{glTK}). More precisely, we define the abstract kernel $K: M \to X'$ by
\begin{align} \label{glK1}
 K(s) := T'\delta_s, 
\end{align} 
where $\delta_s$ is the Dirac functional on $C(M)$ given by $ \left\langle L,\delta_s \right\rangle := L(s) $ for $s \in M$, $L \in C(M)$. An easy computation shows that 
$$ ||K||_\infty = ||T: X \to C(M)||  \quad \text{and} \quad \omega_{X'}(K,\delta) = \omega(T,\delta). $$ 
On the one hand, this implies $K \in C(M,X')$ for $T \in \mathcal K(X,C(M))$, on the other hand we see by the very definition that 
$ |K|_{X',\a} = |T|_\a $. Hence, if $T \in \mathcal Lip_\a(X,C(M))$ is Hölder-continuous of type $\a$, $ 0 < \a \le 1$, then for the kernel $K$ given in (\ref{glK1}) it holds that $K \in C^{\a}(M,X')$ and $||K||_{X',\a} = \Lip_\a(T)$. Moreover, for $x \in X$ and $s \in M$, we have 
$$ (T_K x)(s) = \left\langle x,K(s)\right\rangle = \left\langle x, T'\delta_s \right\rangle = \left\langle Tx,\delta_s \right\rangle = (Tx)(s), $$
which means that the original operator $T$ coincides with the operator $T_K$ generated by the kernel $K$ as given in (\ref{glTK}). Summarizing the above-mentioned facts we arrive at the following well-known statement (cf. \cite[Proposition 5.13.1]{CS90}). 

\begin{prop} \label{propISO}
Let $(M,d)$ be a compact metric space and let $X$ be a Banach space. Then the map $\Phi: \mathcal K(X,C(M)) \to C(M,X')$ defined by 
$$ \Phi(T)(s) = T' \delta_s , \quad s \in M, $$
is a metric isomorphism from $\mathcal K(X,C(M))$ onto $C(M,X')$ as well as a metric isomorphism from the subclass 
$\mathcal Lip_\a(X,C(M))$ of $\mathcal K(X,C(M))$ onto the subclass $C^\a(M,X')$ of $C(M,X')$, for $0<\a\le 1$. 
\end{prop}

Finally, in order to apply our previous results to $C(M)$-valued operators $T_K$ generated by abstract kernels in the sense of (\ref{glTK}) we have to find a link to absolutely convex hulls of precompact sets. To this end, let $S: l_1(M) \to X'$ be the operator defined on the canonical basis $(e_s)_{s \in M}$ of $l_1(M)$ by $Se_s = K(s)$. Furthermore, let $J_\infty$ be the canonical embedding from $C(M)$ into $l_\infty(M)$ and let $\mathcal K_X$ be the canonical metric injection from $X$ into the bidual $X''$. It follows directly from the definitions that 
\begin{align} \label{glS8}
J_\infty T_K = S' \mathcal K_X
\end{align} 
and 
\begin{align} \label{glS7}
e_n(S) = e_n(\aco(\image(K))) \quad \text{and} \quad c_n(S) = c_n(\aco(\image(K))). 
\end{align} 
>From (\ref{glS8}) and the injectivity of the entropy numbers up to the factor two, we obtain
$$ e_n(T_K) \le 2 \, e_n(J_\infty T_K) = 2 \, e_n(S'\mathcal K_X) \le 2 \, e_n(S'). $$ 
The next step is to connect the entropy numbers of $T_K$ with the entropy numbers of the absolutely convex hull of $\image(K)$ by using duality relations. To this end, assume that $X$ is a B-convex Banach space and that the abstract kernel $K:M \to X'$ is continuous. Since also $X'$ is B-convex, the operator $S$ maps into a B-convex Banach space. Furthermore, due to (\ref{glS7}), the operator $S$ is compact. Hence we can use Theorem \ref{endu} to relate the entropy numbers of $S$ with the entropy numbers of $S'$. According to the previous considerations, this relates the entropy of the compact operator $T_K: X \to C(M)$ to the entropy numbers of the absolutely convex hull of $\image(K)$ as a precompact subset of $X'$. In addition, the entropy numbers of $T_K$ can also be related to the Gelfand numbers of $\aco(\image(K))$ by using well known properties of the symmetrized approximation numbers. Indeed, by Theorem 5.3.2 in \cite{CS90} and in view of (\ref{glS8}) and (\ref{glS7}) we have
\begin{align} \label{kolTK}
 d_n(T_K) = t_n(T_K) = t_n(J_\infty T_K) = t_n(S'\mathcal K_X) \le t_n(S') = t_n(S) \le c_n(S) = c_n(\aco(\image(K))) 
\end{align}
and applying Theorem \ref{thC1} to the operator $T_K$ leads to the assertion. Note that we do not need the B-convexity of $X$ in this case. Let us summarize these results in the following lemma.

\begin{lemma} \label{le01}
Let $X$ be a B-convex Banach space. Then for every $0<\a<\infty$ there exists a constant $c=c(\a,X)>0$, such that for the compact operator $T_K:X \to C(M)$ with kernel $K \in C(M,X')$ and all $n \in \N$ we have
$$ \sup_{1 \le k \le n} k^{\a} e_k(T_K) \le c \, \sup_{1 \le k \le n} k^{\a} e_k(\aco(\image(K))). $$
In the case of an arbitrary Banach space $X$, this statement remains true if one replaces $e_k(\aco(\image(K)))$ on the right hand side by the Gelfand numbers $c_k(\aco(\image(K)))$.
\end{lemma} 

As a consequence of Lemma \ref{le01} and the Hardy-type inequalities given in Lemma \ref{LH1} and \ref{LH2} we obtain the following result. The proof is analog to the proof of Theorem \ref{TT03}.  

\begin{coron} \label{cor01} 
Let $X$ be a B-convex Banach space and let $ 0 < r,s < \infty$ and $\a \in \R$. Then there exists a constant $c=c(r,s,\a,X)>0$ such that for the compact operator $T_K: X \to C(M)$ with kernel $K \in C(M,X')$ and all $N \in \N$ it holds
$$ \sum_{n=1}^N \left( \log(n+1) \right)^{\a} \, n^{s/r-1} \, (e_n(T_K))^s \le c \sum_{n=1}^N \left( \log(n+1) \right)^{\a} \, n^{s/r-1} \, (e_n(\aco(\image(K))))^s  $$
and 
$$ \sup_{1 \le n \le N} \left( \log(n+1) \right)^{\a} \, n^{1/r} \, e_n(T_K) \le c \sup_{1 \le n \le N} \left( \log(n+1) \right)^{\a} \, n^{1/r} \, e_n(\aco(\image(K))). $$ 
In particular, we see that 
$$ (e_n(\aco(\image(K))))_n \in l_{r,s,\a} \quad \text{implies} \quad (e_n(T_K))_n \in l_{r,s,\a} $$
for all $0 < r < \infty$, $0 < s \le \infty$ and $\a \in \R$.  
In the case of an arbitrary Banach space $X$, the same holds true for the Gelfand numbers $c_n(\aco(\image(K)))$ instead of $e_n(\aco(\image(K)))$. 
\end{coron} 

Now we are well prepared to give entropy estimates of Lipschitz-continuous operators $T_K : X \to C(M)$ with kernel $K \in C^1(M,X')$ (cf. \cite{CHK88}, \cite[Theorem 5.10.1]{CS90}, \cite{CE01}, \cite{St99}, \cite{St00}, \cite{CHP11}).

\begin{theo} \label{th01}
Let $M$ be a compact metric space with the property that there are constants $\rho,\sigma >0$ and $\gamma \in \R$ such that for all $n \in \N$ it holds
\begin{itemize}
	\item[(a)] $\e_n(M) \le \rho \, n^{-\sigma} \, (\log(n+1))^{-\gamma}$ or 
	\item[(b)] $e_n(M) \le \rho \, n^{-\sigma} \, (\log(n+1))^{-\gamma}$ .
\end{itemize}
Furthermore, let $X$ be a Banach space such that the dual Banach space $X'$ is of type $p$, $1 < p \le 2$. Moreover, let $K \in C^1(M,X')$ be a Lipschitz-continuous kernel and let $T_K:X \to C(M)$ be the corresponding induced operator given by (\ref{glTK}). Then in the case (a) we have the entropy estimate
$$ e_n(T_K) \le c \, n^{-1/p'-\sigma} \, (\log(n+1))^{-\gamma} \, ||K||_{X',1} $$ 
for all $n=1,2,3,\ldots$ and in the case (b) it holds 
\begin{empheq}[left={e_n(T_K) \le c \empheqlbrace}]{alignat*=2}
 & n^{-\sigma} \, (\log(n+1))^{-\gamma} \, ||K||_{X',1} ,                    & \quad & 0<\sigma < 1/p', \\
 & n^{-1/p'} (\log(n+1))^{1/p'-\sigma} (\log(\log(n+3)))^{-\gamma} \, ||K||_{X',1} , & \quad & 1/p' < \sigma < \infty , \\
 & n^{-1/p'} \, (\log(n+1))^{1-\gamma} \, ||K||_{X',1} ,                   & \quad &  \sigma = 1/p', -\infty < \gamma < 1,
\end{empheq}  
for all $n=1,2,3,\ldots$, where $c$ is a constant which may depend on $\rho, \sigma, \gamma, p$ and the type constant of $X$. Furthermore, if $0 < r,s < \infty$, then 
$$ (\e_n(M))_n \in l_{r,s} \quad \text{implies} \quad (e_n(T_K))_n \in l_{q,s} $$  
with $1/q=1/p'+1/r$. Moreover, if $ 0 < r < p'$ and $ 0 < s < \infty$, then 
$$ (e_n(M))_n \in l_{r,s} \quad \text{implies} \quad  (e_n(T_K))_n \in l_{p',s,\alpha} $$ 
with $\alpha = 1/s+1/p'-1/r$ and if $p' < r < \infty$, $0 < s < \infty$ and $\b \in \R$, then 
$$ (e_n(M))_n \in l_{r,s,\b} \quad \text{implies} \quad (e_n(T_K))_n \in l_{r,s,\beta}. $$
\end{theo}
\proof First observe that for all $n \in \N$ it holds that $\e_n(\image(K)) \le ||K||_{X',1} \, \e_n(M)$. In a next step, we use the results of 
\cite[Corollary 5]{St00}, Theorem \ref{TT03}, (\ref{glS5}), (\ref{glS6}),  \cite[Theorem 1.2.]{St04} and Theorem \ref{TT02} to estimate the dyadic entropy numbers of $\aco(\image(K))$. Finally, we carry over these estimates from $e_n(\aco(\image(K)))$ to $e_n(T_K)$ by using Corollary \ref{cor01}. \hfill $\blacksquare$

\begin{rem} \label{rem1} \hfill
\begin{itemize}
\item[(i)] 
Let $(M,d)$ and $X$ be as in Theorem \ref{th01}, but assume that $K: (M,d) \to X' $ is a Hölder-continuous kernel of type $\a$, $ 0 < \a \le 1$. 
Then the estimates in the cases $(a)$ and $(b)$ of Theorem \ref{th01} remain true if we replace the exponents $\sigma$ and $\gamma$ by $\a \sigma$ and $\a \gamma$, respectively. For example, the estimate of the dyadic entropy numbers of $T_K$ in the case $(a)$ reads as
$$ e_n(T_K) \le c \, n^{- 1/p'-\a \sigma} \, (\log(n+1))^{-\a \gamma} \, ||K:(M,d) \to X'||_{X',\a} $$
for all $n=1,2,3,\ldots$. To see this, change the metric on $M$ from $d$ to $d^\a$. Then $K: (M,d^\a) \to X'$ is a Lipschitz-continuous kernel,
$$ || K:(M,d^\a) \to X' ||_{X',1} = ||K:(M,d) \to X'||_{X',\a}, $$
and applying Theorem \ref{th01} with $(M,d^\a)$ yields the assertion due to 
$$ \e_n(M,d^\a) = (\e_n(M,d))^\a. $$
	\item[(ii)] Theorem \ref{th01} remains true for Lipschitz-continuous operators $T \in \mathcal Lip(X,l_\infty(M))$ where $(M,d)$ is a precompact metric space satisfying the entropy condition of Theorem \ref{th01}. This is due to the fact that there exists a compact metric space $(\widehat{M}, \widehat{d})$ and a Lipschitz-continuous operator $S \in \mathcal Lip(X,C(\widehat{M}))$ such that
	$$ \e_n(\widehat{M}, \widehat{d}) \le \e_n(M,d) \le 2 \, \, \e_n(\widehat{M}, \widehat{d}) $$
	and
	$$ \e_n(T: X \to l_\infty(M)) \le \e_n(S:X \to C(\widehat{M})) \le 2 \, \e_n(T: X \to l_\infty(M)) $$
	for all $n \in \N$. 
	\item[(iii)] Consider an operator $T: l_1(M) \to X$, where $(M,d)$ is a precompact metric space and $X$ is a Banach space. Such an operator is said to be Lipschitz-continuous, if its dual $T'$ belongs to $\mathcal Lip(X',l_\infty(M))$ (cf. \cite[Section 3]{CE01}). Theorem \ref{th01} remains true for Lipschitz-continuous operators $T: l_1(M) \to X$, where $X$ is a Banach space of type $p$, $1 < p \le 2$, and $(M,d)$ is a precompact metric space satisfying the entropy condition of Theorem \ref{th01}.  
	\item[(iv)] For later use, we highlight the special case of Theorem \ref{th01} where $X=L_q$ for $1<q<\infty$. 
Then we have in the case $(a)$ 
$$ e_n(T_K) \le c \, n^{-\min\left\{1/q;1/2\right\}-\sigma} \, (\log(n+1))^{-\gamma} \, ||K||_{X',1} $$
for all $n =1,2,3,\ldots$ and in the case $(b)$ we obtain for all $n=1,2,3,\ldots$ the estimates
\begin{align*}
 e_n(T_K) & \le c \, n^{-\sigma} \, (\log(n+1))^{-\gamma} \, ||K||_{X',1} 
\end{align*}
for $\sigma < \min\left\{1/q;1/2\right\}$, 
\begin{align*}
 e_n(T_K) & \le c \, n^{-\min\left\{1/q;1/2\right\}} (\log(n+1))^{\min\left\{1/q;1/2\right\}-\sigma} (\log(\log(n+3)))^{-\gamma} \, ||K||_{X',1}   
\end{align*}
for $\sigma > \min\left\{1/q;1/2\right\} $ and 
\begin{align*}
 e_n(T_K) & \le c \, n^{-\min\left\{1/q;1/2\right\}} \, (\log(n+1))^{-\gamma+1} \, ||K||_{X',1} 
\end{align*}
for $\sigma = \min\left\{1/q; 1/2 \right\}, -\infty < \gamma < 1$.  

\noindent This is due to the fact that $(L_q)'=L_{q'}$ is of optimal type $\min\left\{q';2\right\}$, where $q'$ is given by the Hölder condition $1/q+1/q'=1$.
\end{itemize}
\end{rem}

\subsection{Abstract and integral operators}

Let $M$ be an arbitrary set and let $K:M \to X'$ be a bounded function from $M$ into the dual $X'$ of a Banach space $X$, i.e. 
$$ ||K||_\infty = \sup_{s \in M} ||K(s)||_{X'} < \infty. $$ 
Then $K$ can be considered as an \textit{abstract kernel} which gives rise to an operator $T_K: X \to l_\infty(M)$ from $X$ into the space $l_\infty(M)$ of bounded scalar-valued functions on $M$ by the rule 
\begin{equation} \label{OTK}
(T_K x)(s) := \left\langle x, K(s) \right\rangle, \quad x \in X, s \in M. 
\end{equation}
An easy computation shows that
$$ ||T_K : X \to l_\infty(M) || = ||K||_\infty. $$
Let us define a pseudo-metric on $M$ by  
$$ d(s,t) := ||K(s)-K(t)||_{X'}, \quad s,t \in M.  $$
Next, we introduce the cosets
$$ [t] := \left\{ s \in M : d(s,t) = 0 \right\}, \quad t \in M, $$
and the family of cosets 
$$ \widehat{M} := \left\{ [t] : t \in M \right\}. $$ 
Observe that if $s_1,s_2 \in [s]$ and $t_1,t_2 \in [t]$ then $d(s_1,t_1) = d(s_2,t_2)$. Hence, it makes sense to define a metric $\widehat{d}$ on $\widehat{M}$ by 
$$ \widehat{d}([s],[t]) := d(s,t). $$
It turns out that the entropy numbers of $(M,d)$ and $(\widehat{M},\widehat{d})$ coincide, 
$$ \e_n(M,d) = \e_n(\widehat{M},\widehat{d}). $$

Usually, there is a natural metric $\tilde{d}$ given on $M$. If the entropy numbers of $M$ with respect to the metric $\tilde{d}$ are known and if there is a relationship between $\tilde{d}$ and the pseudo-metric $d$, then it is generally easy to compute the entropy numbers of $M$ with respect to the pseudo-metric $d$ and, therefore, the entropy numbers of $\widehat{M}$ with respect to the metric $\widehat{d}$. As an example, consider the case $M=[0,1]$ with the natural distance $\tilde{d}(s,t) = |s-t|$ and entropy numbers 
$$ \e_n([0,1],\tilde{d}) = (2n)^{-1}. $$
If $d(s,t) \sim \Phi(|s-t|)$ for some continuous strictly increasing function $\Phi:[0,1] \to [0,\infty)$ with $\Phi(0)=0$, then we obtain
$$ \e_n([0,1], d ) \sim \Phi((2n)^{-1}). $$ 

Our aim is to give entropy estimates of the operator $T_K$ mentioned above in (\ref{OTK}). To this end, define the operator $S_K : X \to l_\infty(\widehat{M})$ by
$$ (S_K x)([s]) := (T_Kx)(s), \quad x \in X, [s] \in \widehat{M}, $$
so that $e_n(T_K:X \to l_\infty(M)) = e_n(S_K: X \to l_\infty(\widehat{M}))$, $n \in \N$. Then, for $x \in X$ and $[s],[t] \in \widehat{M}$, it holds
\begin{align*}
 |(S_Kx)([s])-(S_Kx)([t])| & = |(T_Kx)(s)-(T_Kx)(t)| = |\left\langle x, K(s)-K(t) \right\rangle| \\
 &  \le ||x||_X \, ||K(s)-K(t)||_{X'} = ||x||_X \, \widehat{d}([s],[t]),
\end{align*}
which means that $S_K x$ is a continuous function on the metric space $(\widehat{M},\widehat{d})$ and, moreover, that $S_K: X \to l_\infty(\widehat{M})$ is a Lipschitz-continuous operator. Now, given that the metric space $(\widehat{M},\widehat{d})$ is precompact, we can use Theorem \ref{th01} (cf. Remark \ref{rem1} (ii)) to give entropy estimates of the operator $S_K$ and, consequently, of the operator $T_K$.  

In concrete cases, the Banach spaces $X$ and $X'$ are function spaces over some measure space $(\Omega,\mu)$ such that the duality is given by integration with respect to $\mu$,
$$ \left\langle f,g \right\rangle = \int_{\Omega} f(x) g(x) \, \mathrm{d}\mu(x) \quad \text{for } f \in X, g \in X'. $$ 
The kernel $K$ is given as a function $K: M \times \Omega \to \R$ such that $K(s) = K(s,\cdot) \in X'$, $s \in M$. Then the generated operator $T_K:X \to l_\infty(M)$ is given as a kernel integral operator by
$$ (T_K f)(s) = \int_\Omega f(x) \, K(s,x) \, \mathrm{d}\mu(x). $$
To specialize even further, let us now assume that $X=L_p(\Omega, \mu)$ for some $\sigma$-finite measure space $(\Omega,\mu)$ with $1 < p < \infty$. Then $X'=L_{p'}(\Omega,\mu)$ has optimal type $2$, if $1 < p \le 2$ and optimal type $p'$, if $2 < p < \infty$. In this case, the crucial distance on $M$ is given  by
$$ d(s,t) = \left( \int_{\Omega} | K(s,x) - K(t,x) |^{p'} \, \mathrm{d}\mu(x) \right)^{1/p'}. $$

\subsection{Weakly singular integral operators generated by convolution kernels}

In this section we give entropy estimates for weakly singular integral operators 
and weakly singular integral operators of Volterra-type generated by convolution kernels. 
We consider non-negative kernels $K$ on $([0,1] \times [0,1]) \setminus \left\{ (0,0) \right\}$ so that $T_K$ maps a 
function $f$ on $[0,1]$ to the function
$$ (T_K f)(t) = \int_0^1 K(t,x) f(x) \, \mathrm{d}x, \quad t \in [0,1]. $$

In the remainder of this section we distinguish between the cases 
\begin{align*}
 (WS) \qquad & K(t,x) := k(|t-x|),\\
 (VO) \qquad & K(t,x) := \begin{cases} k(t-x), \quad & \text{for } x < t, \\ 0, \quad & \text{for } x \ge t, \end{cases}
\end{align*}
of a weakly singular kernel $(WS)$ and a Volterra-kernel $(VO)$, respectively. 
Here $k: (0,1] \to \R$ is a non-negative, continuous and strictly decreasing function with a singularity at $0$, i.e. 
$$ \lim_{x \to 0} k(x) = \infty. $$
It would actually be enough to assume that $k$ is strictly decreasing only in a neighborhood of $0$. 
Furthermore, we fix $q$ with $1 < q < \infty$ and assume that $k \in L_q[0,1]$.
These are standing assumptions in all of the results to follow. \label{assump}

The following lemma is the key to several examples. 

\begin{lemma} \label{le02}
Let $A \subset [0,1]$ be a subset of the interval $[0,1]$. Denote by $\e_n(A)$ the $n$-th entropy number of $A$ with respect to the absolute value $|\cdot|$. Under the stated assumptions (cf. p.\pageref{assump})
$$ d(s,t) = \left( \int_0^1 |K(s,x)-K(t,x)|^q \, \mathrm{d}x \right)^{1/q} $$
defines a pseudo-metric on $[0,1]$. 
In the case $(WS)$ we have
$$ d(s,t) \le 4^{1/q} \, \left( \int_0^{|s-t|} (k(u))^q \, \mathrm{d}u \right)^{1/q} \quad \text{for } 0 \le s,t \le 1, $$
and 
$$ \e_n(A,d) \le 4^{1/q} \, \left( \int_0^{\e_n(A)} (k(u))^q \, \mathrm{d}u \right)^{1/q} \quad \text{for } n=1,2,3,\ldots . $$
and in the case $(VO)$ it holds 
$$ \left( \int_0^{|s-t|} (k(u))^q \, \mathrm{d}u \right)^{1/q} \le d(s,t) \le 2^{1/q} \, \left( \int_0^{|s-t|} (k(u))^q \, \mathrm{d}u \right)^{1/q} \quad \text{for } 0 \le s,t \le 1, $$
and 
$$ \left( \int_0^{\e_n(A)} (k(u))^q \, \mathrm{d}u \right)^{1/q} \le \e_n(A,d) \le 2^{1/q} \, \left( \int_0^{\e_n(A)} (k(u))^q \, \mathrm{d}u \right)^{1/q} \quad \text{for } n=1,2,3,\ldots. $$
\end{lemma}

\proof In order to estimate the pseudo-metric $d$, let $s>t$. From the inequality 
$$ |a-b|^q \le |a^q-b^q| \quad \text{for } a,b \ge 0 \text{ and } q \ge 1, $$ 
we obtain the estimate
\begin{align*}
(d(s,t))^q & = \int_0^1 |K(s,x)-K(t,x)|^q \, \mathrm{d}x \le \int_0^1 |(K(s,x))^q - (K(t,x))^q| \, \mathrm{d}x \\
& = \int_0^1 |(k(|s-x|))^q - (k(|t-x|))^q| \, \mathrm{d}x  .
\end{align*}
A natural way to proceed is to decompose the integration region. To this end, put 
$$ I:= \int_0^1 |(k(|s-x|))^q - (k(|t-x|))^q| \, \mathrm{d}x = I_1 + I_2 + I_3, $$
where 
\begin{align*}
I_1 & = \int_0^t |(k(|s-x|))^q - (k(|t-x|))^q| \, \mathrm{d}x = \int_0^t \left[ (k(t-x))^q - (k(s-x))^q \right] \, \mathrm{d}x \\
& = \int_0^t (k(u))^q \, \mathrm{d}u - \int_{s-t}^s (k(u))^q \, \mathrm{d}u = \int_0^{s-t} (k(u))^q \, \mathrm{d}u -\int_t^s (k(u))^q \, \mathrm{d}u,
\end{align*}
\begin{align*}
I_2 & = \int_t^s |(k(|s-x|))^q - (k(|t-x|))^q| \, \mathrm{d}x \le \int_t^s \left[ (k(|s-x|))^q + (k(|t-x|))^q \right] \, \mathrm{d}x \\
& = \int_t^s (k(s-x))^q \, \mathrm{d}x + \int_t^s (k(x-t))^q \, \mathrm{d}x = 2 \, \int_0^{s-t} (k(u))^q \, \mathrm{d}u, 
\end{align*}
and 
\begin{align*}
I_3 & = \int_s^1 |(k(|s-x|))^q - (k(|t-x|))^q| \, \mathrm{d}x = \int_s^1 \left[ (k(x-s))^q - (k(x-t))^q \right] \, \mathrm{d}x \\
& = \int_0^{1-s} (k(u))^q \, \mathrm{d}u - \int_{s-t}^{1-t} (k(u))^q \, \mathrm{d}u = \int_0^{s-t} (k(u))^q \, \mathrm{d}u -\int_{1-s}^{1-t} (k(u))^q \, \mathrm{d}u. 
\end{align*}
Consequently, we obtain 
\begin{align*}
(d(s,t))^q & \le I = I_1 + I_2 + I_3 \\
& \le 4 \int_0^{s-t} (k(u))^q \, \mathrm{d}u - \left[ \int_t^s (k(u))^q \, \mathrm{d}u  + \int_{1-s}^{1-t} (k(u))^q \, \mathrm{d}u  \right] \\
& \le 4 \int_0^{s-t} (k(u))^q \, \mathrm{d}u
\end{align*}
for all $0 \le t \le s \le 1$. This implies the desired estimates in the case $(WS)$ of a weakly singular kernel. 
In the case $(VO)$ of a Volterra-kernel 
$$ K(t,x) = \begin{cases} k(t-x), \quad & \text{for } x < t, \\ 0, \quad & \text{for } x \ge t, \end{cases} $$
we get, for $s>t$, that
\begin{align*}
 (d(s,t))^q & = \int_0^1 |K(s,x)-K(t,x)|^q \, \mathrm{d}x  \\ 
& = \int_0^t |k(s-x)-k(t-x)|^q \, \mathrm{d}x + \int_t^s (k(s-x))^q \, \mathrm{d}x .
\end{align*}
Thus, similarly as before, we obtain the estimates
\begin{align*}
\int_t^s (k(s-x))^q \, \mathrm{d}x & \le (d(s,t))^q \\
& \le \int_0^t \left[ (k(t-x))^q-(k(s-x))^q \right] \, \mathrm{d}x + \int_t^s (k(s-x))^q \, \mathrm{d}x \\
& = 2 \int_0^{s-t} (k(u))^q \, \mathrm{d}u - \int_t^s (k(u))^q \, \mathrm{d}u
\end{align*}
and therefore
$$ \int_0^{s-t} (k(u))^q \, \mathrm{d}u \le (d(s,t))^q \le 2 \, \int_0^{s-t} (k(u))^q \, \mathrm{d}u $$
for all $0 \le t \le s \le 1$. This yields the statement and finishes the proof. \hfill $\blacksquare$ 

Now we treat several important examples of weakly singular kernels.

\begin{lemma} \label{le03}
Let $A \subset [0,1]$ be a subset of the interval $[0,1]$. Denote by $\e_n(A)$ the $n$-th entropy number of $A$ with respect to the absolute value $|\cdot|$. Under the stated assumptions (cf. p.\pageref{assump}) the following statements hold: 
\begin{itemize}
	\item[(i)] If the function $k:(0,1] \to \R $ is defined by 
	$$ k(x) = x^{-\tau}, \quad 0<\tau< \frac{1}{q}, $$
	then we have in the case $(WS)$ 
	$$ d(s,t) \preccurlyeq |s-t|^{1/q-\tau} \quad \text{and} \quad \e_n(A,d) \preccurlyeq (\e_n(A))^{1/q-\tau} , $$
	and in the case $(VO)$ it holds
	$$ d(s,t) \sim |s-t|^{1/q-\tau} \quad \text{and} \quad \e_n(A,d) \sim (\e_n(A))^{1/q-\tau}. $$
	\item[(ii)] If the function $k:(0,1] \to \R $ is defined by 
	$$ k(x) = x^{-1/q} (c_0- \ln x)^{-\b}, \quad \frac{1}{q} < \b, $$
	where $c_0$ is a positive constant, then we have in the case $(WS)$ 
	$$ d(s,t) \preccurlyeq (c_0-\ln |s-t|)^{1/q-\b} \quad \text{and} \quad \e_n(A,d) \preccurlyeq (c_0-\ln \e_n(A))^{1/q-\b} , $$
	and in the case $(VO)$ it holds
	$$ d(s,t) \sim (c_0-\ln |s-t|)^{1/q-\b} \quad \text{and} \quad \e_n(A,d) \sim (c_0 - \ln \e_n(A))^{1/q-\b}. $$
\end{itemize}

\end{lemma}

\proof Due to Lemma \ref{le02} it is enough to compute the expression
$$ \left( \int_0^{|s-t|} (k(u))^q \, \mathrm{d}u \right)^{1/q} \quad \text{for } 0 \le s,t \le 1. $$
In the case $(i)$, we have that
$$ \left( \int_0^{|s-t|} (k(u))^q \, \mathrm{d}u \right)^{1/q} = (1- q \tau)^{-1/q} \, |s-t|^{1/q-\tau} \quad \text{for } 0 \le s,t \le 1. $$
Now let us turn to the case $(ii)$. Observe that for any $c_0 > 0$ the function 
$k$ is an element of $L_q[0,1]$ and is strictly decreasing in a neighborhood of $0$. 
The technical assumption that $k$ is strictly decreasing on the whole interval $(0,1]$ 
is ensured if we choose $c_0$ large enough, e.g. $c_0 > \b q$. However, keep in mind 
that the upcoming results are true for any $c_0 >0$. 
For the corresponding integral we get
$$ \left( \int_0^{|s-t|} (k(u))^q \, \mathrm{d}u \right)^{1/q} = (\b q -1 )^{-1/q} \, (c_0-\ln |s-t|)^{1/q-\b} \quad \text{for } 0 \le s,t \le 1. $$
This finishes the proof. \hfill $\blacksquare$

Now we are well prepared to prove the following theorem. Since convolution operators from $L_p[0,1]$ into $C[0,1]$ are, in a sense, closely related to certain diagonal operators from $l_p$ into $l_\infty$, we may expect sharp estimates of entropy numbers of convolution operators only in the case where $2 \le p < \infty$. This is the reason why we restrict our applications to this case. The case $1 < p < 2$ will be treated in another paper. 

\begin{theo} \label{th02}
Under the stated assumptions (cf. p.\pageref{assump}) the following statements hold:
\begin{itemize}
	\item[(i)] If the function $k:(0,1] \to \R $ is defined by 
	$$ k(x) = x^{-\tau}, \quad 2 \le p < \infty, \; 0<\tau< \frac{1}{p'},  $$
	then $T_K$ maps $L_p[0,1]$ into $C[0,1]$ and in the cases $(WS)$ and $(VO)$ the entropy estimate
	$$ e_n(T_K:L_p[0,1] \to C[0,1]) \preccurlyeq n^{\tau - 1} $$
	holds. 
	\item[(ii)] If the function $k:(0,1] \to \R $ is defined by 
	$$ k(x) = x^{-1/p'} (c_0- \ln x)^{-\b}, \quad 2 \le p < \infty, \, \frac{1}{p'} < \b, \; c_0 >0, $$
	then $T_K$ maps $L_p[0,1]$ into $C[0,1]$ and in the cases $(WS)$ and $(VO)$ the following entropy estimates hold: 
	\begin{empheq}[left={e_n(T_K:L_p[0,1] \to C[0,1]) \preccurlyeq \empheqlbrace}]{alignat*=2}
 	  & n^{1/p'-\b},                             & \quad &   1/p' < \b < 1,  \\
 	  & n^{- 1/p} \, (\log(n+1))^{ 1 - \b },     & \quad &   1<\b<\infty,    \\ 
 	  & n^{- 1/p} \, \log(n+1),                  & \quad &   \b = 1. 
	\end{empheq}
\end{itemize}
\end{theo}

\proof The results follow from Lemma \ref{le03} with $A=[0,1]$, $q=p'$ and Remark \ref{rem1} (iv) after Theorem \ref{th01} with $q=p$. 
For the proof of $(i)$ we apply Theorem \ref{th01} $(a)$ with $X=L_p[0,1]$, $\sigma = 1/p' - \tau$ and $ \gamma = 0 $. The proof of $(ii)$ follows from Theorem \ref{th01} $(b)$ with $X=L_p[0,1]$, $\sigma = \b - 1/p'$ and $\gamma = 0$. \hfill $\blacksquare$ 

We can even go a step further and consider more general kernels given by kernel functions 
$$ k(x) = x^{-\tau} \, l(1/x), \quad 0<\tau < \frac{1}{q}, 0 < x \le 1, $$ 
where $l$ is a positive, continuous and slowly varying (in the sense of Karamata) function defined on $[1,\infty)$. Note that such kernel functions always have a singularity at $0$, cf. \cite[Proposition 1.3.6. (v)]{BGT87}. 
Furthermore, \cite[Theorem 1.5.3]{BGT87} tells us that $k(x)$ is up to multiplicative constants equivalent to a decreasing function. This
enables us to apply analogous reasoning as before. 
For more detailed information concerning slowly varying functions consult e.g. Chapter 1 in \cite{BGT87}. 

First we present an analogon of Lemma \ref{le03}.
 
\begin{lemma} \label{leSV}
Under the stated assumptions (cf. p.\pageref{assump}) the following statements hold. 
If the function $k: (0,1] \to \R$ is defined by 
$$ k(x) = x^{-\tau} \, l(1/x), \quad 0 < \tau < 1/q, $$ 
then we have in the case $(WS)$ the estimates
$$ d(s,t) \preccurlyeq |s-t|^{1/q-\tau} \, l(|s-t|^{-1}) \quad \text{and} \quad \e_n([0,1],d) \preccurlyeq n^{\tau-1/q} \, l(2n) $$
and in the case $(VO)$ it holds
$$ d(s,t) \sim |s-t|^{1/q-\tau} \, l(|s-t|^{-1}) \quad \text{and} \quad \e_n([0,1],d) \sim n^{\tau-1/q} \, l(2n). $$
\end{lemma}

\proof Again, it is enough to compute 
$$ \left( \int_0^r (k(u))^q \, \mathrm{d}u \right)^{1/q} \quad \text{for } 0 < r \le 1. $$ 
According to the definition, we have
$$ \int_0^r (k(u))^q \, \mathrm{d}u = \int_0^r u^{-\tau q} \, (l(1/u))^q \, \mathrm{d}u = \int_{1/r}^\infty z^{\tau q -2} \, (l(z))^q \, \mathrm{d}z. $$
Now it follows from \cite[Proposition 1.5.10]{BGT87} that
$$ \int_{1/r}^\infty z^{\tau q -2} \, (l(z))^q \, \mathrm{d}z \sim (1/r)^{\tau q-1} \, (l(1/r))^q. $$
Note that the arising constants depend on $\tau$, $q$ and the function $l$. Hence, 
$$ \left( \int_0^r (k(u))^q \, \mathrm{d}u \right)^{1/q} \sim r^{1/q-\tau} \, l(1/r), $$ 
which yields the assertion. \hfill $\blacksquare$ 

The resulting version of Theorem \ref{th02} then reads as follows. The proof is based on Theorem 6 of \cite{St00}. 

\begin{theo} \label{thSV}
Under the stated assumptions (cf. p.\pageref{assump}) the following statements hold.
If $2 \le p < \infty$ and the function $k: (0,1] \to \R$ is defined by 
$$ k(x) = x^{-\tau} \, l( 1/x ), \quad 0 < \tau < \frac{1}{p'}, $$ 
then in the cases $(WS)$ and $(VO)$ the following entropy estimate holds: For all $\gamma \in \R$ there exists a constant $c=c(\gamma) \ge 1$ such that for all $n \in \N$ we have  
$$ n^{1 - \tau } \, (\log(n+1))^{\gamma} \, e_n(T_K: L_p[0,1] \to C[0,1]) \preccurlyeq c \, \sup_{1 \le k \le n^\b}   (\log(k+1))^\gamma \, l(2k), $$
where $\b = 1 + \frac{p'-1}{1-\tau p'}. $ 

\end{theo}

Finally, let us treat an important example by considering a double-logarithmic term. 

\begin{bsp}
Under the stated assumptions (cf. p.\pageref{assump}) the following statements hold. 
If the function $k: (0,1] \to \R$ is defined by 
$$ k(x) = x^{-\tau} \, (c_0 - \ln x)^{-\b} \, (c_0 + \ln(c_0 - \ln x))^{-\gamma}, \quad 0<\tau \le \frac{1}{q}, \b, \gamma \in \R, $$ 
where $c_0$ is a positive constant, then we have in the case $(WS)$ the estimates
$$ d(s,t) \preccurlyeq f(s,t,\tau,\b,\gamma,c_0,q) \quad \text{and} \quad \e_n([0,1],d) \preccurlyeq g(n,\tau,\beta,\gamma,q) $$
and in the case $(VO)$ we obtain the asymptotic behavior  
$$ d(s,t) \sim f(s,t,\tau,\b,\gamma,c_0,q) \quad \text{and} \quad \e_n([0,1],d) \sim g(n,\tau,\beta,\gamma,q), $$
where $f(s,t,\tau,\b,\gamma,c_0,q) = $
\begin{empheq}[left=\empheqlbrace]{alignat*=2}
 & |s-t|^{1/q-\tau} \, (c_0-\ln |s-t|)^{-\b} \, (c_0 + \ln (c_0 - \ln |s-t|))^{-\gamma} , & \quad & 0<\tau < 1/q, \, \b \in \R, \, \gamma \in \R, \\
 & (c_0-\ln |s-t|)^{1/q-\b} \, (c_0 + \ln (c_0 - \ln |s-t|))^{-\gamma} ,            & \quad & \tau = 1/q, \, 1/q < \b < \infty, \, \gamma \in \R, \\
 & (c_0 + \ln (c_0 - \ln |s-t|))^{1/q-\gamma} ,                                     & \quad & \tau = \b = 1/q, \, 1/q < \gamma < \infty,  
\end{empheq}  
and $g(n,\tau,\beta,\gamma,q) = $
\begin{empheq}[left=\empheqlbrace]{alignat*=2}
 & n^{\tau-1/q} \, (1+\ln n)^{-\b} \, (1 + \ln (1 + \ln n))^{-\gamma}  , & \quad & 0<\tau < 1/q, \, \b \in \R, \, \gamma \in \R , \\
 & (1+\ln n)^{1/q-\b} \, (1 + \ln (1 + \ln n))^{-\gamma} ,               & \quad & \tau = 1/q, \, 1/q < \b < \infty, \, \gamma \in \R, \\
 & (1 + \ln (1 + \ln n))^{1/q-\gamma} ,                                  & \quad & \tau = \b = 1/q, \, 1/q < \gamma < \infty. 
 \end{empheq}  
\end{bsp}

For this example, the resulting version of Theorem \ref{th02} reads as follows.
\begin{prop} \label{th04}
Under the stated assumptions (cf. p.\pageref{assump}) the following statements hold.
If $2 \le p < \infty$ and the function $k: (0,1] \to \R$ is defined by 
$$ k(x) = x^{-\tau} \, (c_0 - \ln x)^{-\b} \, (c_0 + \ln(c_0 - \ln x))^{-\gamma}, \quad 0<\tau \le \frac{1}{p'}, \b, \gamma \in \R, c_0>0, $$ 
then in the cases $(WS)$ and $(VO)$ the following entropy estimates hold: 
$$ e_n(T_K: L_p[0,1] \to C[0,1]) \preccurlyeq f(n,\tau,\beta,\gamma,p), $$ 
where $f(n,\tau,\beta,\gamma,p) = $
\begin{empheq}[left=\empheqlbrace]{alignat=2}
	\label{P1} & n^{\tau-1} \, (\log(n+1))^{-\b} \, (\log \log(n+3))^{-\gamma} ,	& \quad & 0<\tau < 1/p', \, \b \in \R, \, \gamma \in \R,        \\ 
	\label{P2} & n^{1/p' - \b} \, (\log(n+1))^{-\gamma} , 												& \quad & \tau = 1/p', \, 1/p' < \b < 1, \, \gamma \in \R,      \\ 
	\label{P3} & n^{-1/p} \, (\log(n+1))^{1-\b} \, (\log \log(n+3))^{-\gamma},   	& \quad & \tau = 1/p', \, 1 < \b < \infty, \, \gamma \in \R,    \\ 
	\label{P4} & n^{-1/p} \, (\log(n+1))^{1-\gamma}, 														  & \quad & \tau = 1/p', \, \b = 1, \, -\infty < \gamma < 1,      \\ 
	\label{P5} & n^{-1/p} \, (\log(n+1))^{\delta} , 															& \quad & \tau = 1/p', \, \b = 1, \, \gamma \ge 1,  \\ 
	\label{P6} & (\log(n+1))^{1/p'-\gamma} ,											                & \quad & \tau = \b = 1/p', \, 1/p' < \gamma < \infty,        
\end{empheq}
and $\delta$ in $(\ref{P5})$ is an arbitrary positive number.
\end{prop}

\begin{rem} The estimate \eqref{P5} of Proposition \ref{th04} is a consequence of Theorem 1.3. of Creutzig and Steinwart's paper \cite{CrSt02}. To see this, observe that Theorem 1.3. is valid for all $-\infty < \b < 1$ and put $\b = 1- \delta$ for arbitrary $\delta > 0$. 
\end{rem}

Finally, let us deal with the optimality of the results in the Volterra-kernel case $(VO)$. In order to prove lower bounds of the entropy of $T_K: L_p[0,1] \to C[0,1]$ we construct suitable distance nets in $T_K(B_{L_p[0,1]})$. The proofs are inspired by unpublished works of Linde and Lacey in the Hilbert space case $p=2$ (cf. \cite[p. 1807]{Lif10}, \cite[Proposition 38]{Kl12a}). However, it turns out that their techniques also work in a more general framework.

Let us start with (\ref{P1}). Consider the $2^n$ functions 
$$ f_{\e}(x) := \sum_{i=1}^n \e_i \, \mathbbm{1}_{[(i-1)/n , i/n]}(x) , \quad \e = (\e_1,\ldots,\e_n) \in \left\{ -1,1 \right\}^n . $$
Then $f_\e \in B_{L_p[0,1]}$ for every choice of $\e \in \{-1,1\}^n$. Now we estimate the mutual distance of the images of $f_{\e}$ under $T_K$ in the Volterra-kernel case $(VO)$. To this end, let $\e, \tilde{\e} \in \left\{ -1,1 \right\}^n$ with $\e \not= \tilde{\e}$. Let $j$ be the least index such that $\e_j \not= \tilde{\e}_j$. Then 
\begin{align*}
||T_K f_{\e} - T_K f_{\tilde{\e}} ||_\infty & \ge |(T_K f_\e)(j/n) - (T_K f_{\tilde{\e}})(j/n)| = 2 \int_0^{1/n} k(z) \, \mathrm{d}z .
\end{align*} 
Hence, we have found a distance net consisting of $2^n$ elements of $T_K(B_{L_p[0,1]})$ and therefore
$$ e_n(T_K: L_p[0,1] \to C[0,1]) \ge \e_{2^n - 1}(T_K: L_p[0,1] \to C[0,1]) \ge \int_0^{1/n} k(x) \, \mathrm{d}x .$$
Using \cite[Proposition 1.5.10]{BGT87} we compute that  
\begin{align*}
\int_0^{1/n} k(x) \, \mathrm{d}x & = \int_0^{1/n} x^{-\tau} \, (c_0 - \ln x)^{-\b} \, (c_0 + \ln(c_0 - \ln x))^{-\gamma} \, \mathrm{d}x \\
& = \int_n^\infty z^{\tau-2} \, (c_0+\ln z)^{-\b} \, (c_0+\ln(c_0+\ln z))^{-\gamma} \, \mathrm{d}z \\
& \succcurlyeq n^{\tau - 1} \, (c_0+\ln n)^{-\b} \, (c_0+\ln(c_0+\ln n))^{-\gamma} 
\end{align*}
and conclude  
$$ e_n(T_K: L_p[0,1] \to C[0,1]) \succcurlyeq n^{\tau-1} \, (\log(n+1))^{-\b} \, (\log \log(n+3))^{-\gamma}. $$
This shows that estimate (\ref{P1}) is the best possible.

Now let us deal with optimality in the case $\tau = 1/p'$. The idea is to construct a suitable distance net by using the kernel function $k \in L_{p'}[0,1]$. To this end, define functions
\begin{align} \label{fj1} 
 f_j(x) := \frac{1}{\a_m} \, (k(j/m-x))^{p'/p} \, \mathbbm{1}_{[(j-1)/m , \, j/m)}(x), \quad j=1,2,\ldots,m, 
\end{align}
where 
$$ \a_m = \left( \int_0^{1/m} (k(x))^{p'} \, \mathrm{d}x \right)^{1/p} . $$
Then $f_j \in B_{L_p[0,1]}$ for every $j=1,2,\ldots,m$. Furthermore, for $1 \le i < j \le m$, we obtain 
\begin{align*}
 ||T_K f_i - T_K f_j ||_\infty & \ge |(T_K f_i)(i/m) - (T_K f_j)(i/m)| = \frac{1}{\a_m} \int_{0}^{1/m} (k(z))^{p'} \, \mathrm{d}z = \a_m^{p-1}.
\end{align*}
Hence, 
$$ \e_{m-1}(T_K: L_p[0,1] \to C[0,1]) \ge \frac{1}{2} \, \a_m^{p-1}. $$
Now let $\tau = 1/p'$, $\b > 1/p'$ and $\gamma \in \R$. Using \cite[Proposition 1.5.10]{BGT87} we compute that
\begin{align*}
\a_m^{p-1} & = \left( \int_0^{1/m} (k(x))^{p'} \, \mathrm{d}x \right)^{1/p'} = \left( \int_0^{1/m} x^{-1} \, (c_0 - \ln x)^{-\b p'} \, (c_0 + \ln(c_0 - \ln x))^{-\gamma p'} \, \mathrm{d}x \right)^{1/p'} \\
& = \left( \int_{c_0 + \ln m}^\infty z^{-\b p'} (c_0 + \ln z)^{- \gamma p'} \, \mathrm{d}z \right)^{1/p'} \succcurlyeq (c_0 + \ln m)^{1/p'-\b} \, (c_0 + \ln (c_0 + \ln m))^{-\gamma}	 .
\end{align*}
Consequently, putting $m=2^{n-1}+1$, we get 
$$ e_n(T_K: L_p[0,1] \to C[0,1]) \succcurlyeq n^{1/p'-\b} \, (\log(n+1))^{-\gamma} $$ 
This shows that estimate (\ref{P2}) is the best possible. Moreover, we see that in the critical case $\tau = 1/p'$ and $\b = 1$ the estimate 
$$ e_n(T_K: L_p[0,1] \to C[0,1]) \succcurlyeq n^{-1/p} (\log(n+1))^{-\gamma} $$ 
holds for $\gamma \in \R$, cf. (\ref{P4}). 

Now let us deal with the case $\tau = \b = 1/p'$ and $\gamma > 1/p'$. Here we have 
\begin{align*}
 \a_m^{p-1}  = (\gamma p' - 1)^{-1/p'} \, (c_0 + \ln(c_0 + \ln m))^{1/p'-\gamma} 
\end{align*} 
and therefore
$$ e_n(T_K:L_p[0,1] \to C[0,1]) \succcurlyeq (\log(n+1))^{1/p'-\gamma}. $$
This shows that estimate (\ref{P6}) is the best possible.

Finally, we show that estimate (\ref{P3}) is the best possible. To see this, we consider suitable means of the functions $f_j$ defined in (\ref{fj1}). For $J \subset \left\{ 1,2,\ldots,m \right\}$ define
$$ f_J(x) := |J|^{-1/p} \sum_{j \in J} f_j(x) .$$
Then $f_J \in B_{L_p[0,1]}$ for every choice of $J \subset \left\{ 1,2,\ldots,m \right\}$. Let $m>1$ be a square number and define
$$ \Phi_m := \left\{ f_J : J \subset \left\{1,2,\ldots,m \right\} \text{ with } |J| = \sqrt{m} \right\}. $$
Then 
$$ \log_2 |\Phi_m| = \log_2 \binom{m}{\sqrt{m}} \ge \log_2\left( \left(\frac{m}{\sqrt{m}}\right)^{\sqrt{m}} \right) = \frac{1}{2} \sqrt{m} \, \log_2(m) .$$
Let $f_J, f_L \in \Phi_m$ with $J \not= L$ and let $i$ be the least element in the symmetric difference $(J \cup L) \setminus (J \cap L)$. Then, for $\tau=1/p'$, $\b > 1/p'$ and $\gamma \in \R$, we have
\begin{align*}
|| T_K f_J - T_K f_L ||_\infty & \ge | (T_K f_J)(i/m) - (T_K f_L)(i/m) | \\
& = (\sqrt{m})^{-1/p} \, \a_m^{-1} \int_0^{1/m} (k(z))^{p'} \, \mathrm{d}z \\
& = (\sqrt{m})^{-1/p} \, \a_m^{p-1} \\
& \succcurlyeq (\sqrt{m})^{-1/p} \, (c_0 + \ln m)^{1/p'-\b} \, (c_0 + \ln (c_0 + \ln m))^{-\gamma} \\
& \succcurlyeq \left( \sqrt{m} \, \log_2(m) \right)^{-1/p} \, \left( \log\left(\sqrt{m} \log_2(m)+ 1 \right)  \right)^{1-\b} \times \\
& \quad \times \left( \log\log\left( \sqrt{m} \, \log_2(m) +3 \right) \right)^{-\gamma} 
\end{align*} 
Hence, we have found at least $2^{\frac{1}{2} \sqrt{m} \, \log_2(m)}$ elements in $T_K(B_{L_p[0,1]})$ with mutual distance (up to some constant) at least  
$$ \left( \sqrt{m} \, \log_2(m) \right)^{-1/p} \, \left( \log\left(\sqrt{m} \log_2(m)+ 1 \right)  \right)^{1-\b} \, \left( \log\log\left( \sqrt{m} \, \log_2(m) +3 \right) \right)^{-\gamma}, $$
where $m>1$ is an arbitrary square number. Therefore,
$$ e_n(T_K: L_p[0,1] \to C[0,1] ) \succcurlyeq n^{-1/p} \, (\log(n+1))^{1-\b} \, (\log\log(n+3))^{-\gamma} $$
and this shows that estimate (\ref{P3}) is the best possible.

\subsection{Weakly singular integral operators from $L_2[0,1]$ in $C[0,1]$ and $L_q[0,1]$ generated by convolution kernels}

Given a precompact subset $A$ of a Banach space $X$ of type $p$, we do not have exact entropy estimates of $\aco(A)$ in the critical case that
$$ e_n(A) \preccurlyeq n^{-1/p'} \, (\log(n+1))^{\b} \quad \text{with } \b \ge 1. $$ 
In contrast to that, in the Hilbert space case we have such estimates and very recent developments show, that they are asymptotically optimal. Consequently, some of our estimates of $e_n(T_K)$ from the previous section can be refined in the Hilbert space case. This fact and the general importance of the Hilbert space case motivates this section. We start with recalling one of the main results of \cite{CHP11} which gives a complete overview about entropy and Gelfand numbers of absolutely convex hulls in the Hilbert space case (see also \cite[Theorem 1]{Kl12b}). 

\begin{theo} \label{enHil}
Let $(s_k)$ stand for the Gelfand numbers $(c_k)$ or for the dyadic entropy numbers $(e_k)$. Let $A \subset H$ be a precompact subset of a Hilbert space $H$. If $0<r<\infty$ and $\b \in \R$ then there exists a constant $c=c(r,\b,A) > 0$ such that the following inequality holds in the respective cases:
\begin{enumerate}
\item[(i)] If $0<r<2$ and $\b \in \R$ then 
	\begin{align*}
		\sup_{1 \le k < \infty} & (\log \log(k+3))^\b \, (\log(k+1))^{1/r-1/2} \, k^{1/2} \, s_k(\aco(A)) \\
		& \le c \left( 1 + \sup_{1 \le k < \infty} (\log(k+1))^\b \, k^{1/r} \, e_k(A) \right) .	
	\end{align*}
\item[(ii)] If $r=2$ then  
	\begin{align*}
		\sup_{1 \le k \le n} (\log(k+1))^{\b-1} \, k^{1/2} \, s_k(\aco(A)) \le c \left( 1 + \sup_{1 \le k \le n} (\log(k+1))^\b \, k^{1/2} \, e_k(A) \right) 
	\end{align*}
	for $-\infty < \b < 1$ and $n \in N$, and 
	\begin{align*}
		\sup_{1 \le k \le n} (\log \log(k+3))^{-1} \, k^{1/2} \, s_k(\aco(A)) \le c \left( 1 + \sup_{1 \le k \le n} \log(k+1) \, k^{1/2} \, e_k(A) \right) 
	\end{align*}
	for $\b=1$ and $n \in \N$, and
	\begin{align*}
		\sup_{1 \le k < \infty} (\log \log(k+3))^{\b-1} \, k^{1/2} \, s_k(\aco(A)) \le c \left( 1 + \sup_{1 \le k < \infty} (\log(k+1))^\b \, k^{1/2} \, e_k(A) \right) 
	\end{align*}
	for $1<\b<\infty$. 	
\item[(iii)] If $2<r<\infty$ and $\b \in \R$ then the expressions
	\begin{align*}
		\sup_{1 \le k \le n} (\log(k+1))^\b \, k^{1/r} \, s_k(\aco(A)) \quad \text{and} \quad 1+\sup_{1 \le k \le n} (\log(k+1))^\b \, k^{1/r} \, e_k(A)
	\end{align*}
	are asymptotically equivalent.
\end{enumerate}
\end{theo}

It is well known that the results given in Theorem \ref{enHil} $(i)$ and $(iii)$ are the best possible ones (cf. \cite{CHP11}). In the critical case $(ii)$, where
$$ e_n(A) \preccurlyeq n^{-1/2} \, (\log(n+1))^{\b} \quad \text{with } \b \in \R, $$
it is known for almost ten years that the result for $ - \infty < \b < 1$ is asymptotically optimal. This goes back to Gao \cite{G01} for $\b=0$. Creutzig and Steinwart \cite{CrSt02} extended Gao's ideas to $ - \infty < \b < 1$ and to $B$-convex Banach spaces. Very recent results of Gao \cite{G12} show that also the results for $\b=1$ and $\b>1$ are asymptotically optimal. This is the subject of the next theorem. 

\begin{theo} Let $H$ be an infinite dimensional separable Hilbert space. 
\begin{enumerate}
	\item[(i)] {\em \cite{G01}} There exists a subset $A \subset H$ of $H$ and positive constants $c_1,c_2$ such that 
		$$ \log N(A,\e) \le c_1 \, \e^{-2} \quad \text{for all } \e > 0 $$
		and 
		$$ \log N(\aco(A),\e) \ge c_2 \, \e^{-2} \, |\log \e|^2 \quad \text{for all } 0< \e < c_2. $$
		This means that there exists a subset $A \subset H$ of $H$ with
		$$ e_n(A) \preccurlyeq n^{-1/2}  \quad \text{and} \quad e_n(\aco(A)) \succcurlyeq n^{-1/2} \, \log(n+1) .$$		
	\item[(ii)] {\em \cite{CrSt02}} For $\b < 1$, there exists a subset $A \subset H$ of $H$ with
		$$ e_n(A) \preccurlyeq n^{-1/2} \, (\log(n+1))^{-\b}  \quad \text{and} \quad e_n(\aco(A)) \succcurlyeq n^{-1/2} \, (\log(n+1))^{1-\b} .$$		
	\item[(iii)] {\em \cite{G12}} There exists a subset $A \subset H$ of $H$ and positive constants $c_1,c_2$ such that 
		$$ \log N(A,\e) \le c_1 \, \e^{-2} \, |\log \e|^{-2} \quad \text{for all } 0<\e<1/2 $$
		and 
		$$ \log N(\aco(A),\e) \ge c_2 \, \e^{-2} \, (\log |\log \e|)^2 \quad \text{for all } 0< \e < 2^{-6}. $$
		This means that there exists a subset $A \subset H$ of $H$ with
		$$ e_n(A) \preccurlyeq n^{-1/2} \, (\log(n+1))^{-1}  \quad \text{and} \quad e_n(\aco(A)) \succcurlyeq n^{-1/2} \, \log \log(n+3) .$$
	\item[(iv)] {\em \cite{G12}} For $\b > 1$, there exists a subset $A \subset H$ of $H$ and positive constants $c_1,c_2$ such that 
		$$ \log N(A,\e) \le c_1 \, \e^{-2} \, |\log \e|^{-2\b} \quad \text{for all } 0<\e<1/2 $$
		and 
		$$ \log N(\aco(A),\e) \ge c_2 \, \e^{-2} \, (\log |\log \e|)^{2-2\b} \quad \text{for all } 0< \e < 2^{-6}. $$
		This means that, for $\b > 1$, there exists a subset $A \subset H$ of $H$ with
		$$ e_n(A) \preccurlyeq n^{-1/2} \, (\log(n+1))^{-\b}  \quad \text{and} \quad e_n(\aco(A)) \succcurlyeq n^{-1/2} \, (\log \log(n+3))^{1-\b} .$$
\end{enumerate}
\end{theo}

In (\ref{kolTK}) we related the Kolomogorov numbers $d_n(T_K)$ of the operator $T_K$ to the Gelfand numbers $c_n(\aco(\image(K)))$ of the absolutely convex hull of $\image(K)$. This relationship leads to fruitful results in the Hilbert space case. Indeed, we can use Theorem \ref{enHil} to relate the entropy numbers of $\image(K)$ to both $e_n(\aco(\image(K)))$ and $c_n(\aco(\image(K)))$. Hence, in the Hilbert space case we can estimate not only the dyadic entropy numbers $e_n(T_K)$ but also the Kolmogorov numbers $d_n(T_K)$ of the operator $T_K$. In the Hilbert space setting, the resulting version of Proposition \ref{th04} reads as follows.

\begin{theo} \label{enTkH}
Let $(s_n)$ stand for the Kolmogorov numbers $(d_n)$ or for the dyadic entropy numbers $(e_n)$. Under the stated assumptions (cf. p.\pageref{assump}) the following statements hold. If the function $k: (0,1] \to \R$ is defined by 
$$ k(x) = x^{-\tau} \, (c_0 - \ln x)^{-\b} \, (c_0 + \ln(c_0 - \ln x))^{-\gamma}, \quad 0<\tau \le \frac{1}{2}, \b, \gamma \in \R, c_0>0, $$ 
then in the cases $(WS)$ and $(VO)$ the following estimates hold: 
$$ s_n(T_K:L_2[0,1] \to C[0,1]) \preccurlyeq  f(n,\tau,\beta,\gamma) , $$
where $f(n,\tau,\beta,\gamma) =$
\begin{empheq}[left=\empheqlbrace]{alignat=2}
	\label{G1} & n^{\tau-1} \, (\log(n+1))^{-\b} \, (\log \log(n+3))^{-\gamma} ,	& \quad & 0<\tau < 1/2, \, \b \in \R, \, \gamma \in \R,           \\ 
	\label{G2} & n^{1/2 - \b} \, (\log(n+1))^{-\gamma} , 												 	& \quad & \tau = 1/2,   \, 1/2 < \b < 1, \, \gamma \in \R,        \\ 
	\label{G3} & n^{-1/2} \, (\log(n+1))^{1-\b} \, (\log \log(n+3))^{-\gamma},   	& \quad & \tau = 1/2,   \, 1 < \b < \infty, \, \gamma \in \R,     \\ 
	\label{G4} & n^{-1/2} \, (\log(n+1))^{1-\gamma}, 														  & \quad & \tau = 1/2,   \, \b = 1, \, -\infty < \gamma < 1,       \\ 
	\label{G5} & n^{-1/2} \, \log \log(n+3) , 																		& \quad & \tau = 1/2,   \, \b = 1, \, \gamma = 1,                 \\ 
	\label{G6} & n^{-1/2} \, (\log \log(n+3))^{1-\gamma} ,											  & \quad & \tau = 1/2,   \, \b = 1, \, 1 < \gamma < \infty,        \\ 
  \label{G7} & (\log(n+1))^{1/2-\gamma} , 																			& \quad & \tau = \beta = 1/2, \, 1/2 < \gamma < \infty.   	
\end{empheq}
 
\end{theo}

According to Theorem \ref{enTkH}, the  behavior of entropy numbers as well as Kolmogorov numbers of the operator $T_K: L_2[0,1] \to C[0,1]$ differs significantly between the cases $0<\tau<1/2$, $\tau=1/2, \, \b > 1/2$ and $\tau = \beta = 1/2$. Furthermore, we see that for fixed $\tau = 1/2$ a sudden jump occurs if the parameter $\b$ crosses the point $\b=1$. In addition, for fixed $\tau=1/2$ and $\b=1$, we recognize a sudden jump if the parameter $\gamma$ crosses the point $\gamma = 1$. 

We already know from the previous subsection that the entropy estimates given in (\ref{G1}), (\ref{G2}), (\ref{G3}) and (\ref{G7}) are the best possible. In the critical case (\ref{G4}), Lifshits \cite[Theorem 3.2]{Lif10} proved that, for $\gamma = 0$,
$$ e_n(T_K: L_2[0,1] \to C[0,1]) \preccurlyeq n^{-1/2} . $$
Hence, in the critical case, our general approach using absolutely convex hulls does not lead to a sharp upper estimate. We do not know whether the upper estimates given in $(\ref{G5})$ and $(\ref{G6})$ are the best possible. Furthermore, we would like to point out that Linde \cite{Lin08} proved the lower estimate 
$$ e_n(T_K:L_2[0,1] \to C[0,1]) \succcurlyeq n^{-1/2} \, (\log(n+1))^{1/2-\b} $$
in the case where $\tau = 1/2$, $\b > 1/2$ and $\gamma = 0$.

Without proof we remark that the estimates (\ref{G1}), (\ref{G2}), (\ref{G3}) and (\ref{G7}) of the Kolmogorov numbers of $T_K$ are also optimal. This can be derived from the optimality of the entropy estimates in this cases by using Theorem \ref{thC1} in combination with a trick given in \cite[p. 106]{C85}.

In contrast to Theorem \ref{enTkH} we now study entropy and Kolmogorov numbers of convolution operators from $L_2[0,1]$ into $L_q[0,1]$ for $1 \le q < \infty$. It turns out that the asymptotic behavior of those numbers significantly changes in the critical cases. This demonstrates the difficulties of estimating entropy and Kolmogorov numbers of convolution operators. 

Let us start with recalling the $l$-norm of an operator $T: X \to Y$ (or absolutely $\gamma$-summing norm in \cite{LP74}). Let $l_2^n$ be the $n$-dimensional Euclidean space and $S: l_2^n \to Y$ an operator, then the \textit{$l$-norm of $S$} is defined by
$$ l(S) := \left( \int_{R^n} ||Sx||^2 \, \mathrm{d}\gamma_n(x) \right) ^{1/2}, $$
where $\gamma_n$ is the canonical Gaussian probability measure of $\R^n$. For an operator $T: X \to Y$ we define 
$$ l(T) := \sup \left\{ l(TA) : ||A: l_2^n \to X || \le 1, n \in \N \right\} . $$
If $A: X_0 \to X$ and $B : Y \to Y_0$ are operators acting between Banach spaces, then $l$ has the ideal property (cf. \cite{LP74})
$$ l(BTA) \le ||B|| \, l(T) \, ||A|| . $$
Furthermore, we need a refined version of a Sudakov-type inequality. The following theorem is due to Pajor and Tomczak-Jaegermann. 

\begin{theo} \label{suda} { \em \cite{PT86} }
There is a constant $c \ge 1$ such that for all operators $T: X \to H$ from a Banach space $X$ into a Hilbert space $H$ and all $n \in N$,
\begin{align} \label{sda01}
 n^{1/2} \, c_n(T) \le c \, l(T') . 
\end{align}
\end{theo} 

\noindent By Gordon \cite{Go88} we know that $c \le \sqrt{2}$. 

For our purposes, we need an additional version of Pajor and Tomczak-Jaegermann's inequality (see also \cite[Lemma A]{CE03}). In order to formulate it we introduce the approximation numbers with respect to the $l$-norm. For an operator $T: X \to Y$ acting between Banach spaces $X$ and $Y$ the \textit{approximation numbers with respect to the $l$-norm} are defined by
$$ a_n(T ; l) := \inf \{ l(T-A) : A \in \mathcal L(X,Y) \text{ with } \rank A < n \}, \quad n=1,2,\ldots. $$
Analogously, we define $a_n(T ; \Pi_q)$ as the \textit{approximation numbers with respect to the absolutely $q$-summing norm} (cf. \cite{P87}).  

The approximation numbers with respect to the $l$-norm were used for some time in functional analysis with different notations (cf. e.g. \cite[Theorem 9.1]{Pi89}).
They also play a role in probability theory since they describe the approximability of Gaussian processes by finite sums (cf. e.g. \cite{LL99}).

\begin{lemma} \label{suda2}
\begin{enumerate}
	\item[(i)] For an operator $T: H \to X$ from a Hilbert space $H$ into a Banach space $X$ with $l(T)<\infty$ we have the inequality 
		$$ k^{1/2} d_{k+n-1}(T) \le \sqrt{2} \, a_n(T; l) \quad \text{for } k,n \in \N. $$
	\item[(ii)] Let $1 \le q < \infty$, then for all absolutely $q$-summing operators $T: H \to X$ from a Hilbert space $H$ into a Banach space $X$ we have the inequality
	$$ k^{1/2} d_{k+n-1}(T) \le \sqrt{2q} \, a_n(T; \Pi_q) \quad \text{for } k,n \in \N. $$
\end{enumerate}
\end{lemma}

\proof Since $l(T)<\infty$, the operator $T: H \to X$ is compact and, therefore, we have that $d_n(T) = c_n(T')$ (cf. \cite[11.7.7]{P78}, \cite[Proposition 2.5.6]{CS90}). Furthermore, it holds that $T'' = \mathcal K_X T \mathcal K_H^{-1}$, where $\mathcal K_Z$ is the canonical metric injection from a Banach space $Z$ into its bidual $Z''$. We conclude that $l(T'') \le l(T)$ and taking (\ref{sda01}) into account gives
$$ n^{1/2} d_n(T) = n^{1/2} c_n(T') \le \sqrt{2} \, l(T'') \le \sqrt{2} \, l(T) , \quad n \in \N.  $$
Now let $A: H \to X$ be an operator with $\rank(A) < n$. Due to the additivity and rank property of the Kolmogorov numbers, we get
$$ d_{k+n-1}(T) \le d_k(T-A) + d_n(A) = d_k(T-A) $$
for all $k,n \in \N$. Hence, we obtain
$$ k^{1/2} d_{k+n-1}(T) \le k^{1/2} d_k(T-A) \le \sqrt{2} \, l(T-A) $$
for all $k,n \in \N$ and all operators $A: H \to X$ with $\rank(A) < n$. This yields the assertion
$$ k^{1/2} d_{k+n-1}(T) \le \sqrt{2} \, a_n(T;l) \quad \text{for } k,n \in \N. $$
Now let us deal with the proof of $(ii)$. By Linde and Pietsch \cite{LP74} we have for an absolutely $q$-summing operator $T: H \to X$ the estimate
$$ l(T) \le b_q \, \Pi_q(T), $$ 
where 
$$ b_q = \max \left\{ 1 ; 2^{1/2} \frac{\Gamma(\frac{q+1}{2})}{\Gamma(\frac{1}{2})}   \right\} \le \sqrt{q}, \quad 1 \le q < \infty. $$
Combining this estimate with $(i)$, we get the desired assertion. \hfill $\blacksquare$

Now we are well prepared to prove the following theorem. 

\begin{theo} \label{rieli}
	Let $k : (0,1] \to \R$ be a kernel function as stated on page \pageref{assump} with $k \in L_2[0,1]$. Then for the weakly singular integral operator $T_K : L_2[0,1] \to L_q[0,1]$, $1 \le q < \infty$, given by
	$$ (T_K f)(t) = \int_0^1 k(|t-x|) \, f(x) \, \mathrm{d}x $$
the inequality 
$$ n^{1/2} d_n(T_K: L_2[0,1] \to L_q[0,1]) \le c \, \sqrt{q} \, \left( \int_0^{1/n} (k(u))^2 \, \mathrm{d}u \right)^{1/2}, \quad n \in \N, $$
holds true, where $c \ge 1$ is an absolute constant.   
\end{theo}

\proof By Lemma \ref{le02} we have with 
$$ d(s,t) = \left( \int_0^1 |K(s,x)-K(t,x)|^2 \, \mathrm{d}x \right)^{1/2} $$
the estimate 
$$ \e_n([0,1],d) \le 2 \, \left( \int_0^{1/{2n}} (k(u))^2 \, \mathrm{d}u \right)^{1/2} \quad \text{for } n=1,2,3,\ldots . $$
Using \cite[Theorem 5.6.1]{CS90} (see also \cite{RS96}) we get 
$$ a_{n+1}(T_K : L_2[0,1] \to C[0,1]) \le \e_n([0,1],d) \le 2 \, \left( \int_0^{1/{2n}} (k(u))^2 \, \mathrm{d}u \right)^{1/2}. $$
Moreover, for the identity operator $I: C[0,1] \to L_q[0,1]$ we have that $\Pi_q(I)=1$ (cf. \cite[1.3.8]{P87}). Using the inequality
$$ a_{n+1}(T_K : L_2[0,1] \to L_q[0,1] ; \Pi_q) \le \Pi_q(I: C[0,1] \to L_q[0,1]) \, a_{n+1}(T_K : L_2[0,1] \to C[0,1]) $$
we arrive at
$$ a_{n+1}(T_K : L_2[0,1] \to L_q[0,1] ; \Pi_q) \le 2 \left( \int_0^{1/{2n}} (k(u))^2 \right)^{1/2} . $$
Combining this estimate with Lemma \ref{suda2} $(ii)$ we finally obtain
$$ k^{1/2} d_{k+n}(T_K : L_2[0,1] \to L_q[0,1]) \le 2  \sqrt{2q} \, \left( \int_0^{1/{2n}} (k(u))^2 \right)^{1/2} $$
for $k,n \in \N$. Putting $k=n$ and $k=n-1$, respectively, we get with a new absolute constant $c \ge 1$ the desired estimate
$$ n^{1/2} d_n(T_K : L_2[0,1] \to L_q[0,1]) \le c \, \sqrt{q} \left( \int_0^{1/{n}} (k(u))^2 \right)^{1/2} $$
for $n=1,2,\ldots$. \hfill $\blacksquare$

Now we give the corresponding result to Theorem \ref{enTkH} for weakly singular convolution operators from $L_2[0,1]$ into $L_q[0,1]$. 

\begin{theo} \label{enTkH2}
Let $(s_n)$ stand for the Kolmogorov numbers $(d_n)$ or for the dyadic entropy numbers $(e_n)$. Under the stated assumptions (cf. p.\pageref{assump}) the following statements hold. If the function $k: (0,1] \to \R$ is defined by 
$$ k(x) = x^{-\tau} \, (c_0 - \ln x)^{-\b} \, (c_0 + \ln(c_0 - \ln x))^{-\gamma}, \quad 0<\tau \le \frac{1}{2}, \b, \gamma \in \R, c_0>0, $$ 
then in the cases $(WS)$ and $(VO)$ the following estimates hold for all $1 \le q < \infty$: 
$$ s_n(T_K:L_2[0,1] \to L_q[0,1]) \le c(\b,\gamma) \, \sqrt{q} \,  f(n,\tau,\beta,\gamma) , $$
where $f(n,\tau,\beta,\gamma) =$
\begin{empheq}[left=\empheqlbrace]{alignat=2}
	\label{J1} & n^{\tau-1} \, (\log(n+1))^{-\b} \, (\log \log(n+3))^{-\gamma} ,	& \quad & 0<\tau < 1/2, \, \b \in \R, \, \gamma \in \R,           \\ 
	\label{J2} & n^{-1/2} \, (\log(n+1))^{1/2-\b} \, (\log \log(n+3))^{-\gamma}, 	& \quad & \tau = 1/2,   \, 1/2 < \b < \infty, \, \gamma \in \R,   \\ 
	\label{J3} & n^{-1/2} \, (\log \log(n+3))^{1/2-\gamma},                      	& \quad &  \tau = \beta = 1/2, \, 1/2 < \gamma < \infty.     
\end{empheq}
\end{theo}

\proof Using Theorem \ref{rieli} we obtain the desired estimates for the Kolmogorov numbers of $T_K$. By applying Theorem \ref{thC1} (see also \cite[Theorem 1.3]{CKP99}) we get the same asymptotic estimates also for the dyadic entropy numbers of $T_K$. \hfill $\blacksquare$

If we compare the estimates of entropy and Kolmogorov numbers of weakly singular convolution operators from $L_2[0,1]$ into $C[0,1]$ given in Theorem \ref{enTkH} with those of Theorem \ref{enTkH2} 
for the same convolution operator considered from $L_2[0,1]$ into $L_q[0,1]$
we observe a significant difference in the critical case $\tau = 1/2$, $\b > 1/2$, $\gamma \in \R$  and in the super-critical case $\tau = \beta = 1/2, \, \gamma > 1/2$. 
In the critical case the difference is on the logarithmic scale, in the super-critical case it is even on the polynomial scale.
In particular, we see that the estimates given in Theorem \ref{enTkH} are not the limiting case $q \to \infty$ of the estimates in Theorem \ref{enTkH2}.

\subsection{Riemann-Liouville operator}
In this section we deal with entropy and Kolmogorov numbers of the \textit{famous Riemann-Liouville operator}
$$ (R_\a f)(t) := \frac{1}{\Gamma(\a)} \int_0^t (t-x)^{\a-1} \, f(x) \, \mathrm{d}x $$
for $0 \le t \le 1$ and $\a > 0$. Singular numbers (= approximation numbers) of these operators between Hilbert spaces have been extensively studied by many authors in the literature (cf. e.g. \cite{Bu07,DM97,Do93,Do95,FM86,M01}). Our results extend and complement results in the literature, especially of Lomakina and Stepanov \cite{LS06}, Li and Linde \cite{LL99} and Linde \cite{Lin04}. We start with recalling the fact that the classical Riemann-Liouville operator satisfies the semigroup property
$$ R_\a (R_\b f) = R_{\a+\b} f , \quad \a,\b > 0 . $$
Furthermore, we need the following result for a general Volterra integration operator.

\begin{lemma} \label{RL01}
Let $k: (0,1] \to \R$ be a kernel function defined as on page \pageref{assump} with $k \in L_1[0,1]$. Then for the Volterra-operator 
$$ (T_Kf)(t) = \int_0^t k(t-x) f(x) \, \mathrm{d}x $$
we have that
$$ T_K : L_p[0,1] \to L_p[0,1] , \quad 1 \le p \le \infty, $$
and for all $0 \le \delta \le 1$ the estimate 
$$ || (T_Kf)(\cdot+\delta) - (T_Kf)(\cdot) | L_p[0,1] || \le 2 \, ||f | L_p[0,1]||  \int_0^{\delta} k(s) \, \mathrm{d}s    $$
holds true. 
\end{lemma}

\proof Using the triangle inequality in $L_p$ we see that
\begin{align*}
& || (T_Kf)(\cdot+\delta) - (T_Kf)(\cdot)||_p \\ 
& = \left( \int\limits_0^{1-\delta} \left| (T_Kf)(t+\delta) - (T_Kf)(t)  \right|^p \mathrm{d}t \right)^{1/p} \\
& = \left( \int\limits_0^{1-\delta} \, \left| \int\limits_0^{t+\delta} k(t+\delta-x) f(x) \, \mathrm{d}x - \int\limits_0^t k(t-x) f(x) \, \mathrm{d}x \right|^p \mathrm{d}t \right)^{1/p} \\
& = \left( \int\limits_0^{1-\delta} \, \left| \int\limits_0^t \left[ k(t+\delta-x) - k(t-x) \right] f(x) \, \mathrm{d}x + \int\limits_t^{t+\delta} k(t+\delta-x) f(x) \, \mathrm{d}x \right|^p \mathrm{d}t \right)^{1/p} \\
& \le \left( \int\limits_0^{1-\delta} \, \left| \int\limits_0^t \left[ k(t+\delta-x) - k(t-x) \right] f(x) \, \mathrm{d}x \right|^p \mathrm{d}t \right)^{1/p} + \left( \int\limits_0^{1-\delta} \, \left| \int\limits_t^{t+\delta} k(t+\delta-x) f(x) \, \mathrm{d}x \right|^p \mathrm{d}t \right)^{1/p} \\
& =: I_1 + I_2 .
\end{align*} 
Now we estimate the first integral $I_1$. Applying Minkowski's integral inequality (cf. \cite[Theorem 202]{HLP88}) we obtain
\begin{align*}
I_1 & = \left( \int\limits_0^{1-\delta} \, \left| \int\limits_0^t \left[ k(s+\delta) - k(s) \right] f(t-s) \, \mathrm{d}s \right|^p \mathrm{d}t \right)^{1/p} \\
& \le \left( \int\limits_0^{1-\delta} \, \left( \int\limits_0^{1-\delta} \left| k(s+\delta) - k(s) \right| |f(t-s)| \, \mathbbm{1}_{[0,t]}(s) \, \mathrm{d}s \right)^p \mathrm{d}t \right)^{1/p} \\
& \le \int\limits_0^{1-\delta} \left( \int\limits_{s}^{1-\delta} \left| k(s+\delta)-k(s) \right|^p |f(t-s)|^p \, \mathrm{d}t  \right)^{1/p} \mathrm{d}s  \\
& = \int\limits_0^{1-\delta} \left| k(s+\delta)-k(s) \right| \left( \int\limits_{s}^{1-\delta} |f(t-s)|^p \, \mathrm{d}t \right)^{1/p} \mathrm{d}s \\
& \le ||f||_p \, \int\limits_0^{1-\delta} \left| k(s+\delta)-k(s) \right| \, \mathrm{d}s .
\end{align*}
Since $k$ is non-negative and decreasing we can continue with
\begin{align*}
\int_0^{1-\delta} \left| k(s+\delta)-k(s) \right| \, \mathrm{d}s & = \int_0^{1-\delta} k(s) \, \mathrm{d}s - \int_0^{1-\delta} k(s+\delta) \, \mathrm{d}s \\
& = \int_0^{\delta} k(s) \, \mathrm{d}s - \int_0^{\delta} k(1-s) \, \mathrm{d}s \le \int_0^{\delta} k(s) \, \mathrm{d}s. 
\end{align*}
Hence, we have that
$$ I_1 \le ||f||_p \, \int_0^{\delta} k(s) \, \mathrm{d}s . $$
Finally we deal with the second integral $I_2$. Using again Minkowski's integral inequality we see that
\begin{align*}
I_2 & \le \left( \int\limits_0^{1-\delta} \, \left( \int\limits_0^{\delta} k(\delta-s) \, |f(s+t)| \, \mathrm{d}s \right)^p \mathrm{d}t \right)^{1/p}  \\
& \le \int_0^{\delta} \left( \int_{0}^{1-\delta} (k(\delta-s))^p \, |f(s+t)|^p \, \mathrm{d}t  \right)^{1/p} \mathrm{d}s \\
& \le ||f||_p \, \int_0^{\delta} k(\delta-s) \, \mathrm{d}s = ||f||_p \, \int_{0}^\delta k(s) \, \mathrm{d}s .
\end{align*}
The claim follows. \hfill $\blacksquare$ 

Now we are able to give upper entropy estimates of the classical Riemann-Lioville operator.

\begin{prop} \label{RL02} Let $1 \le p,q \le \infty$ and $\max\{ 1/p-1/q ; 0 \} < \a \le 1$. Then
$$ e_n(R_\a : L_p[0,1] \to L_q[0,1]) \preccurlyeq n^{-\a} . $$
\end{prop}

\proof Applying Lemma \ref{RL01} with the kernel function $k(x) = \frac{1}{\Gamma(\a)} x^{\a-1}$ gives
$$ || (R_\a f)(\cdot+\delta) - (R_\a f)(\cdot) | L_p[0,1] || \le \frac{2}{\Gamma(\a)} \, ||f | L_p[0,1]||  \int_0^{\delta} s^{\a-1} \, \mathrm{d}s = \frac{2}{\a \Gamma(\a)} \delta^\a \, ||f | L_p[0,1]||. $$
This implies that the image $R_\a(L_p[0,1])$ belongs to the Besov space $B_{p,\infty}^\a [0,1]$. Now we factorize the Riemann-Liouville operator $R_\a: L_p[0,1] \to L_q[0,1]$ as $R_\a = I_{p,q} \, R_\a^B$, where $R_\a^B : L_p[0,1] \to B_{p,\infty}^\a [0,1]$ is the operator $R_\a$ with values in the Besov space and $I_{p,q}: B_{p,\infty}^\a [0,1] \to L_q[0,1]$ is the natural embedding. From \cite{C81b} (see also \cite{ET96}, \cite{Kö86}) we know that 
$$  e_n(I_{p,q}: B_{p,\infty}^\a [0,1] \to L_q[0,1]) \sim n^{-\a} , \quad 1 \le p,q \le \infty.  $$
Consequently, we get the desired upper estimate
$$ e_n(R_\a : L_p[0,1] \to L_q[0,1]) \le || R_\a^B|| \,  e_n(I_{p,q}) \preccurlyeq n^{-\a} . $$ \hfill $\blacksquare$

In a next step, using Proposition \ref{RL02} we give a complete overview about the entropy behavior of the classical Riemann-Liouville operator. Our results extend those of \cite[Theorem 2.2]{LS06} from integers to positive real numbers. 

\begin{theo} \label{RL03}
Let $1 \le p,q \le \infty$ and $ \max\{1/p-1/q ; 0\} < \a < \infty$. Then 
$$ e_n(R_\a : L_p[0,1] \to L_q[0,1]) \sim n^{-\a} . $$
\end{theo}

\proof The upper estimate is a consequence of Proposition \ref{RL02} and \cite[Theorem 2.2]{LS06} in combination with the semigroup property of the Riemann-Liouville operator. First, let $1< p,q <\infty$. In view of Proposition \ref{RL02} and \cite[Theorem 2.2]{LS06} it is enough to consider the case $\a > 1$ and $\a \notin \N$. We write $\a = \left\lfloor \a \right\rfloor + a$, $0<a<1$. Due to \cite[Theorem 2.2]{LS06} we have 
$$ e_n( R_{\left\lfloor \a \right\rfloor} : L_p[0,1] \to L_q[0,1] ) \preccurlyeq n^{-\left\lfloor \a \right\rfloor} $$
for $1 < p,q < \infty$ and Proposition \ref{RL02} gives
$$ e_n(R_a : L_p[0,1] \to L_p[0,1] ) \preccurlyeq n^{-a} .$$
>From the multiplicativity of entropy numbers we conclude that
$$ e_{2n-1}(R_\a) = e_{2n-1}(R_{\left\lfloor \a \right\rfloor} R_a) \le e_n(R_{\left\lfloor \a \right\rfloor}) \, e_n(R_a) \preccurlyeq n^{-\left\lfloor \a \right\rfloor} \, n^{-a} = n^{-\a} $$
for all $n=1,2,\ldots$. This yields the desired upper estimate for $1 < p,q < \infty$. It remains to deal with the cases $p \in \{ 1,\infty \}$, $1 \le q \le \infty$ and $1 \le p \le \infty$, $q \in \{ 1,\infty \}$. Let, for example, $p=1$ and $1<q<\infty$. For $\a > 1$, we write $\a = a + b$ with $b > 1$ and $0<a<1$. Choose $r > 1$ such that $1-1/r < a$ and consider $R_a : L_1[0,1] \to L_r[0,1]$ and $R_b : L_r[0,1] \to L_q[0,1]$. Due to Proposition \ref{RL02} we have that $ e_n(R_a : L_1[0,1] \to L_r[0,1]) \preccurlyeq n^{-a} $ and we have already proved that $ e_n(R_b : L_r[0,1] \to L_q[0,1]) \preccurlyeq n^{-b}  $.
Hence, 
$$ e_{2n-1}(R_\a : L_1[0,1] \to L_q[0,1]) \preccurlyeq n^{-a} \, n^{-b} = n^{-\a} . $$
The other cases can be treated similarly. 

In order to prove the lower bound we use \cite[Theorem 2.2]{LS06} and the upper estimate of $e_n(R_\a)$ we have already proved. First, let $1 < p,q < \infty$. For $\a > \max \{ 1/p-1/q ; 0 \}$ and $\a \notin \N$ we choose $0<a<1$ such that $\left\lfloor \a \right\rfloor + 1 = \a + a $. Then \cite[Theorem 2.2]{LS06} gives
$$ e_n(R_{\left\lfloor \a \right\rfloor + 1} : L_p[0,1] \to L_q[0,1]) \succcurlyeq n^{-\left\lfloor \a \right\rfloor - 1} .$$
Consequently, applying Proposition \ref{RL02} we obtain
\begin{align*}
	n^{-\left\lfloor \a \right\rfloor - 1} & \preccurlyeq e_{2n-1}(R_{\left\lfloor \a \right\rfloor + 1} : L_p[0,1] \to L_q[0,1]) \\
	& \le e_n(R_a : L_p[0,1] \to L_p[0,1]) \, e_n(R_\a : L_p[0,1] \to L_q[0,1]) \\
	& \preccurlyeq n^{-a} \, e_n(R_\a : L_p[0,1] \to L_q[0,1]) 
\end{align*}
and conclude that
$$ e_n(R_\a : L_p[0,1] \to L_q[0,1]) \succcurlyeq n^{-\left\lfloor \a \right\rfloor - 1 + a} = n^{-\a}. $$
The remaining cases $p \in \{ 1,\infty \}$, $1 \le q \le \infty$ and $1 \le p \le \infty$, $q \in \{ 1,\infty \}$ can be treated similarly. \hfill $\blacksquare$ 

\begin{rem}
	Similar statements as in Theorem \ref{RL03} can be obtained for a multiplicative $s$-number sequence $(s_n)$ if the $s$-numbers of the embedding 
	$$ I_{p,q}: B_{p,\infty}^\a[0,1] \to L_q[0,1] ,\quad 1 \le p,q \le \infty, \quad \a > \max \{1/p-1/q;0\}, $$
	are known. 
\end{rem}

Finally, we deal with entropy and Kolmogorov numbers of weakly singular integral operators $T_K: L_p[0,1] \to C(A)$, $2 \le p < \infty$, where $A \subset [0,1]$ is a compact subset of the interval $[0,1]$. The inequalities of the following theorem are modifications of the inequalities of \cite[Theorem 6]{St00} and \cite{CE03}.

\begin{theo} \label{RL04}
Let $2 \le p < \infty$ and $A \subset [0,1]$ be a compact subset of the interval $[0,1]$. Furthermore, let $k : (0,1] \to \R$ be a kernel function as stated on page \pageref{assump} with $k \in L_{p'}[0,1]$. Then for the weakly singular integral operator $T_K : L_p[0,1] \to C(A)$,
	$$ (T_K f)(t) = \int_0^1 k(|t-x|) \, f(x) \, \mathrm{d}x , \quad t \in A, $$
the following inequalities hold true: 
\begin{enumerate}
	\item[(i)] 
	  \begin{align*} \sup_{1 \le k \le n} k^{\rho + 1/p} \, (\log(k+1))^\gamma \, e_k(T_K: L_p[0,1] \to C(A)) \\
	  \le c(\rho,\gamma,p,K) \sup_{1 \le k \le n^{1+\frac{1}{p \rho}}} k^{\rho} \, (\log(k+1))^\gamma \, \e_k(A,d)  
	 \end{align*}
		for $\rho >0$, $\gamma \in \R$ and $n \in \N$, where $ \e_n(A,d) \le 4^{1/p'} \left( \int_0^{\e_n(A)} (k(u))^{p'} \, \mathrm{d}u \right)^{1/p'} $. In the case $p=2$ the inequality holds also for the Kolmogorov numbers $d_n(T_K: L_2[0,1] \to C(A))$. 
	\item[(ii)] $$ n^{1/2} d_n(T_K: L_2[0,1] \to C(A)) \le c(K) \left( 1 + \sum_{k=1}^n k^{-1/2} \, e_k(A,d) \right)  $$
	and 
	$$ k^{1/2} d_{k+n}(T_K : L_2[0,1] \to C(A)) \le c \left( (\log(n+1))^{1/2} \, \e_n(A,d) + \sum_{j=n+1}^\infty \frac{\e_j(A,d)}{j \, (\log(j+1))^{1/2}}  \right) $$
	for $k,n \in \N$, where $ \e_n(A,d) \le 2 \left( \int_0^{\e_n(A)} (k(u))^{2} \, \mathrm{d}u \right)^{1/2} $ and $e_n(A,d) = \e_{2^{n-1}}(A,d)$.	In particular, we have that
	$$ 2^{n/2} \, d_{2^n}(T_K : L_2[0,1] \to C(A)) \le c \left( n^{1/2} \, e_n(A,d) + \sum_{j=n}^\infty j^{-1/2} \, e_j(A,d) \right) $$
	for $n \in \N$.
\end{enumerate}
\end{theo}

\proof From \cite[Theorem 6]{St00} and Corollary \ref{cor01} we get the desired estimate $(i)$ (see also Lemma \ref{le02}). In this context we remark that Theorem $6$ of \cite{St00} is valid for all $\gamma \in \R$. Furthermore, from \ref{kolTK} we know that $d_n(T_K) \le c_n(\aco(\image(K)))$, where the kernel $K: A \to L_{p'}[0,1]$ maps $A$ into $L_{p'}[0,1]$. In the Hilbert space case $p=2$, we use the inequalities of \cite[Theorem B]{CE03} to obtain the desired result $(ii)$ (see also \cite[p. 7]{CHP11}). \hfill $\blacksquare$

Theorem \ref{RL04} is the key to several estimates of entropy and Kolmogorov numbers of weakly singular integral operators. In the following, we complement results for the classical Riemann-Liouville operator given in \cite{Lin04}.

\begin{prop} \label{RL05}
	Let $A \subset [0,1]$ be a compact subset of the interval $[0,1]$ with $\e_n(A) \preccurlyeq n^{-\delta} \, (\log(n+1))^\theta$ for $\delta \ge 1$ and $\theta \in \R$ (note that since $A \subset [0,1]$ we necessarily have $\theta \le 0$ for $\delta = 1$). Then for the classical Riemann-Liouville operator $R_\a : L_p[0,1] \to C(A)$, $2 \le p < \infty$, $\a > 1/p$, we have the estimate
	$$ e_n(R_\a : L_p[0,1] \to C(A)) \preccurlyeq n^{-1/p - \delta(\a - 1/p)} \, (\log(n+1))^{\theta(\a-1/p)} $$
	and in the case $p=2$, $\a > 1/2$, we also obtain
	$$ d_n(R_\a : L_2[0,1] \to C(A)) \preccurlyeq n^{-1/2 - \delta(\a - 1/2)} \, (\log(n+1))^{\theta(\a-1/2)}. $$
\end{prop}

\proof For the kernel function $k(x) = \frac{1}{\Gamma(\a)} x^{\a-1}$, $\a > 1/p$, we obtain 
\begin{align*} 
\e_n(A,d) & \le 4^{1/p'} \left( \int_0^{\e_n(A)} (k(u))^{p'} \, \mathrm{d}u \right)^{1/p'} = \frac{1}{\Gamma(\a)} \left( \frac{4}{p'(\a-1)+1} \right)^{1/p'} \, (\e_n(A))^{\a-1/p} \\
 & \preccurlyeq n^{-\delta(\a-1/p)} \, (\log(n+1))^{\theta(\a-1/p)} . 
\end{align*}
Now we apply Theorem \ref{RL04} $(i)$ with $\rho = \delta(\a - 1/p)$ and $\gamma = - \theta (\a - 1/p) $ to obtain the desired estimate 
$$ e_n(R_\a : L_p[0,1] \to C(A)) \preccurlyeq n^{-1/p - \delta(\a - 1/p)} \, (\log(n+1))^{\theta(\a-1/p)} . $$
In the case $p=2$ we get the corresponding estimate also for the Kolmogorov numbers of $R_{\a}: L_2[0,1] \to C(A)$. \hfill $\blacksquare$

In order to illustrate the generality of the inequalities given in Theorem \ref{RL04} we prove an analogon of Theorem \ref{enTkH} for the limiting case $\tau = 1/2$ and $1/2 < \b < \infty$.  

\begin{prop} \label{RL06}
Let $A \subset [0,1]$ be a compact subset and $T_K : L_2[0,1] \to C(A)$ the weakly singular integral operator generated by the kernel function 
	$$ k(x) = x^{-1/2} \, (c_0 - \ln x)^{-\b} , \quad 1/2 < \beta < \infty, c_0 >0 . $$
Furthermore, let $(s_n)$ stand for the Kolmogorov numbers $(d_n)$ or for the dyadic entropy numbers $(e_n)$. Then the following statements hold true:
\begin{enumerate}
	\item[(i)] If $\e_n(A) \sim n^{-\delta}$ for $\delta \ge 1$ then 
	\begin{empheq}[left={s_n(T_K : L_2[0,1] \to C(A)) \preccurlyeq \empheqlbrace}]{alignat=2}
	 & n^{1/2-\b}  ,											& \quad & 1/2 < \b <1,         \nonumber  \\ 
	 & n^{-1/2} \, \log(n+1), 						& \quad & \b = 1,   \nonumber \\ 
	 & n^{-1/2} \, (\log(n+1))^{1-\b},    & \quad & 1 < \b < \infty. \nonumber    
	\end{empheq}
	\item[(ii)] If $\e_n(A) \sim e^{-n^\delta}$ for $\delta > 0$ then
	$$ s_n(T_K : L_2[0,1] \to C(A)) \preccurlyeq n^{-1/2 - \delta(\b-1/2)}. $$
\end{enumerate}
\end{prop}

\proof Both results follow from Theorem \ref{RL04}. At first, we compute that
$$ \e_n(A,d) \le 2 \left( \int_0^{\e_n(A)} (k(u))^{2} \, \mathrm{d}u \right)^{1/2} \preccurlyeq (c_0 - \ln \e_n(A))^{1/2-\b} . $$ 
Consequently, in the case $(i)$ we obtain $ \e_n(A,d) \preccurlyeq (\log(n+1))^{1/2-\b} $ and $e_n(A,d) \preccurlyeq n^{1/2-\b}$. Applying the inequalities of Theorem \ref{RL04} $(ii)$ leads to the desired estimates for the Kolmogorov numbers of $T_K$. For example, if $1/2 < \b < 1$ then 
$$ \sum_{k=1}^n k^{-1/2} \, e_k(A,d) \preccurlyeq \sum_{k=1}^n k^{-\b} \preccurlyeq n^{1-\b} $$
and this yields 
$$ n^{1/2} \, d_n(T_K : L_2[0,1] \to C(A)) \preccurlyeq n^{1-\b} . $$
The other cases can be treated similarly. Using Theorem \ref{thC1} (see also \cite[Theorem 1.3]{CKP99}) we obtain the same asymptotic estimates also for the dyadic entropy numbers of $T_K$. In the case $(ii)$ we get $\e_n(A,d) \preccurlyeq n^{\delta(1/2-\b)}$. Consequently, applying Theorem \ref{RL04} $(ii)$ with $\rho = -\delta(1/2-\b)$ and $\gamma=0$ yields
$$ s_n(T_K : L_2[0,1] \to C(A)) \preccurlyeq n^{-1/2+\delta(1/2-\b)} $$
and finishes the proof. \hfill $\blacksquare$ 

\begin{rem}
	The upper estimates given in Proposition \ref{RL06} $(i)$ are the same as for the whole interval $A=[0,1]$, cf. Theorem \ref{enTkH}. The parameter $\delta$ only affects the constant. 
\end{rem}

We can also prove an analogon of Theorem \ref{rieli}. 

\begin{theo} \label{rieli2}
	Let $k : (0,1] \to \R$ be a kernel function as stated on page \pageref{assump} with $k \in L_2[0,1]$. Furthermore, let $A \subset [0,1]$ be a compact subset of the interval $[0,1]$ and $\mu$ a Hausdorff measure on $A$. Then for the weakly singular integral operator $T_K : L_2[0,1] \to L_q[A,\mu]$, $1 \le q < \infty$, given by
	$$ (T_K f)(t) = \int_0^1 k(|t-x|) \, f(x) \, \mathrm{d}x $$
the estimate
$$ n^{1/2} d_{2n}(T_K: L_2[0,1] \to L_q[A,\mu]) \le c \, (\mu(A))^{1/q} \, \sqrt{q} \, \left( \int_0^{\e_n(A)} (k(u))^2 \, \mathrm{d}u \right)^{1/2}, \quad n \in \N, $$
holds true, where $c \ge 1$ is an absolute constant. In particular, for $k(x) = \frac{1}{\Gamma(\a)} x^{\a-1}$ , $\a > 1/2$, we obtain that
$$ d_{2n}(T_K : L_2[0,1] \to L_q[A,\mu]) \preccurlyeq n^{-1/2} \, (\e_n(A))^{\a-1/2}. $$ 
\end{theo}

\proof The proof is analog to that of Theorem \ref{rieli}. This time we use that
$$ a_{n+1}(L_2[0,1] \to C(A)) \le \e_n(A,d) \le 2 \left( \int_0^{\e_n(A)} (k(u))^2 \, \mathrm{d}u \right)^{1/2} $$
and
$$ \Pi_q(I: C(A) \to L_q[A,\mu]) \le (\mu(A))^{1/q}, \quad 1 \le q < \infty. $$
In case that $k(x) = \frac{1}{\Gamma(\a)} x^{\a-1}$, $\a > 1/2$, we compute that 
$$ \left( \int_0^{\e_n(A)} (k(u))^2 \, \mathrm{d}u \right)^{1/2} = \frac{1}{\Gamma(\a)} \left( \frac{(\e_n(A))^{2(\a-1)+1}}{2(\a-1)+1} \right)^{1/2}. $$
The assertion follows. \hfill $\blacksquare$

\subsection{Concluding remarks}

Finally, we discuss the relation between Gelfand numbers and Gelfand widths of operators (cf. \cite[p. 336]{P07}). 
In the literature, the use of those notions sometimes lead to confusion. The following comments are designed to clarify this subject.

Let $A \subset Y$ be a bounded subset of a Banach space $Y$. Then the \textit{$n$-th Gelfand width} of the subset $A$ is defined by
$$ c^n(A) := \inf \left\{ \sup \left\{ ||y|| : y \in A \cap F \right\}  : F \text{ subspace of } Y \text{ with } \codim(F) < n \right\}. $$ 
So $c^n(A)$ measures the slices $A \cap F$ of $A$ with subspaces $F$ of codimension at most $n-1$ so that the maximal norm of the slice is as small as possible. 
The $n$-th Gelfand width of $A$ can also be characterized as follows
$$ c^n(A) = \inf \left\{ \sup \left\{ ||y|| : y \in A, \left\langle y, y_k' \right\rangle = 0 \text{ for all } 1 \le k < n  \right\}  : y_1',\ldots,y_{n-1}' \in Y' \right\} . $$
Now we mention two concepts to carry over this notion to linear bounded operators $T: X \to Y$ between Banach spaces $X$ and $Y$. The first one was introduced by Triebel \cite{Tr70}. He defined the \textit{$n$-th Gelfand width} of an operator $T \in \mathcal L(X,Y)$ by
$$ c^n(T) := c^n(T(B_X)), $$
which can be written as
\begin{align} \label{CR01}
 c^n(T) = \inf \left\{ \sup \left\{ ||Tx|| : x \in B_X, \left\langle x, T'y_k' \right\rangle = 0 \text{ for all } 1 \le k < n  \right\}  : y_1',\ldots,y_{n-1}' \in Y' \right\}  . 
\end{align}
The second concept goes back to Pietsch \cite{P74}. He defined the \textit{$n$-th Gefand number} of an operator $T \in \mathcal L(X,Y)$ as in the introduction by 
$$ c_n(T) := \inf \left\{ \sup \left\{ ||Tx|| : x \in B_X \cap E \right\}  : E \text{ subspace of } X \text{ with } \codim(E) < n \right\}. $$ 
This can be written as 
\begin{align} \label{CR02}
 c_n(T) = \inf \left\{ \sup \left\{ ||Tx|| : x \in B_X, \left\langle x, x_k' \right\rangle = 0 \text{ for all } 1 \le k < n  \right\}  : x_1',\ldots,x_{n-1}' \in X' \right\}  .
\end{align}
The main difference between both notions is that intersection is done in the target space $Y$ in the first case while it is done in the domain space $X$ in the second case. 
Furthermore, we remark that the Gelfand numbers of operators $c_n(T)$ are injective but not surjective whereas the Gelfand widths of operators $c^n(T)$ are injective and surjective. 

In view of (\ref{CR01}) and (\ref{CR02}), we immediately check that the Gelfand widths of an operator $T$ in the sense of Triebel and the Gelfand numbers of an operator $T$ in the sense of Pietsch coincide in the case that the dual operator $T'$ is surjective. Otherwise, one only has that 
$$ c_n(T) \le c^n(T). $$ 
In general, this inequality is strict and the difference between $c_n(T)$ and $c^n(T)$ can be quite large. To give an example, we consider the operators 
$$ T_n : l_1 \stackrel{Q_n}{\longrightarrow} l_2^n \stackrel{I_n}{\longrightarrow} l_\infty^n, \quad n=1,2,3,\ldots, $$
as given in Pietsch's monograph \cite[p. 336]{P07}, where $Q_n: l_1 \to l_2^n$ is a metric surjection and $I_n: l_2^n \to l_\infty^n$ is the identity map. Then for the $n$-th Gelfand widths and the $n$-th Gelfand numbers of the operators $T_{2n}$ the estimates
$$ c^n(T_{2n}) \ge \frac{1}{\sqrt{2}} \quad \text{and} \quad c_n(T_{2n}) = a_n(T_{2n}) = d_n(T_{2n}) \sim n^{-1/2} $$
hold. To see this, we use the following arguments:  

Since the Gelfand widths $c^n$ are surjective and the dual of $I_{n}: l_2^{n} \to l_\infty^{n}$ is a surjective operator, we obtain by a result of Ste\v{c}kin (cf. \cite[11.11.8]{P78}) and from \cite[11.5.3, 11.7.4]{P78} that
\begin{align*}
 c^n(T_{2n}) & = c^n(I_{2n}Q_{2_n}) = c^n(I_{2n}:l_2^{2n} \to l_\infty^{2n}) = c_n(I_{2n}:l_2^{2n} \to l_\infty^{2n}) \\
 & = a_n(I_{2n}:l_2^{2n} \to l_\infty^{2n}) =  \sqrt{\frac{n+1}{2n}} \ge \frac{1}{\sqrt{2}} .
\end{align*}
On the other hand, the operator $T_{2n}$ maps $l_1$ into an $l_\infty$-space, and therefore the Gelfand, Kolmogorov and approximation numbers of $T_{2n}$ coincide (cf. \cite[11.5.3, 11.6.3]{P78}). Hence, using a duality argument (cf. \cite[11.7.6]{P78}) and a result of Kashin \cite{Ka77} and Garnaev/Gluskin \cite{GG84} (see \cite{CP88} for a generalization) we see that
$$ c_n(T_{2n}) = a_n(T_{2n}) = d_n(T_{2n}) = d_n(I_{2n}:l_2^{2n} \to l_\infty^{2n}) \sim n^{-1/2}. $$
Both estimates show that the Gelfand widths of an operator can be much larger than the Gelfand numbers of an operator and even larger than the approximation numbers. Consequently, we see that the Gelfand widths cannot be $s$-numbers in the sense of Pietsch. By the way, the difference between $c_n(T)$ and $c^n(T)$ was stressed for the first time in \cite[p. 30]{Pin85}. 

\begin{rem} \hfill
\begin{itemize}

\item[(i)] We should emphasize that if we consider the $n$-th Gelfand width and the $n$-th Gelfand number of the operator $T_n$, then we get equality. Indeed, using the same arguments as above we find that
$$ c^n(T_n) = c^n(I_n:l_2^n \to l_\infty^n) = c_n(I_n:l_2^n \to l_\infty^n)=a_n(I_n:l_2^n \to l_\infty^n)= n^{-1/2} $$ 
and
$$ c_n(T_n) = a_n(T_n) = d_n(T_n) = d_n(I_n:l_2^n \to l_\infty^n) = n^{-1/2}. $$
Hence, $c^n(T_n) = c_n(T_n)$ in this case. The equality
$$ d_n(I_n:l_2^n \to l_\infty^n) = n^{-1/2} $$
results from the following general fact: If $(s_n)_n$ is an $s$-number sequence in the sense of Pietsch \cite{P87} with the property that for the identity operator $I: X \to X$ on an $n$-dimensional Banach space $X$ the $n$-th $s$-number is one, $s_n(I) = 1$, then for any invertible operator $T: X \to Y$ between arbitrary $n$-dimensional Banach spaces the formula $s_n(T) ||T^{-1}|| = 1$ is true. In particular, this gives for the Kolmogorov numbers the equality
$$ d_n(I_n:l_2^n \to l_\infty^n) \, ||I_n^{-1}: l_\infty^n \to l_2^n || = d_n(I_n:l_2^n \to l_\infty^n) \, n^{1/2} = 1. $$

\item[(ii)] If we consider the absolutely convex hull $\aco(A)$ of a bounded subset $A \subset X$ of a Banach space $X$ then, in general, we have (cf. (\ref{glaco}) and (\ref{glS4}))
$$ c_n(\aco(A)) = c_n(T_A: l_1(A) \to X) \le c^n(T_A : l_1(A) \to X) = c^n(\aco(A)). $$
In the case where the dual operator $(T_A)'$ of $T_A:l_1(A) \to X$ is surjective we get equality,
$$ c_n(\aco(A)) = c^n(\aco(A)). $$ 
For example, this applies for finite subsets $A \subset X$ consisting of linearly independent elements. Finally, we point out that one can characterize the $n$-th Gelfand width of $\aco(A)$ by 
$$ c^n(\aco(A)) = \frac{1}{2} \inf\left\{ \diam( \aco(A) \cap E ) : E \text{ subspace of } X \text{ with } \codim(E) < n \right\} , $$
where $\diam( \aco(A) \cap E ) = \sup \left\{ ||x-y|| : x,y \in \aco(A) \cap E \right\} $ denotes the diameter of $\aco(A) \cap E$. In particular, if $X$ is a finite dimensional Banach space with $\dim X = N$ then we get
$$ c^n(\aco(A)) = \frac{1}{2} \inf\left\{ \diam( \aco(A) \cap E ) : E \text{ subspace of } X \text{ with } \dim(E) > N-n \right\}. $$
This means that $c^n(\aco(A))$ measures the minimal diameter of $m$-dimensional slices of $\aco(A)$, $m > N-n$.  

\end{itemize}
\end{rem}


\begin{thebibliography}{ABCDEF} 

\bibitem[AMST04]{AMST04} S. Artstein, V. Milman, S. Szarek, N. Tomczak-Jaegermann, On convexified packing and entropy duality. Geom. Funct. Anal. \textbf{14} (2004), 1134--1141

\bibitem[BP90]{BP90} K. Ball, A. Pajor, The entropy of convex bodies with few extreme points. In: Geometry of Banach spaces (Strobl, 1989), London Math. Soc. Lecture Note Ser. \textbf{158}, 25--32, Cambridge Univ. Press, Cambridge, 1990

\bibitem[B32]{B32} S. Banach, Th\'eorie des Op\'erations Lin\'eaires, Warsaw, 1932

\bibitem[Bu07]{Bu07} P. Burman, Sharp bounds for singular values of fractional integral operators. J. Math. Anal. Appl. \textbf{327} (1) (2007), 251--256

\bibitem[BGT87]{BGT87} N.H. Bingham, C.M. Goldie, J.L. Teugels, Regular Variation. Cambridge Univ. Press, Cambridge, 1987

\bibitem[BPST89]{BPST89} J. Bourgain, A. Pajor, S.J. Szarek, N. Tomczak-Jaegermann, On the duality problem for entropy numbers of operators. Geometric aspects of functional analysis (1987-88), 50--63, Lecture Notes in Math. \textbf{1376} (1989), Springer Berlin

\bibitem[C81a]{C81a} B. Carl, Entropy numbers, s-numbers, and eigenvalue problems. J. Funct. Anal. \textbf{41} (1981), 290--306

\bibitem[C81b]{C81b} B. Carl, Entropy numbers of embedding maps between Besov spaces with an application to eigenvalue problems. Proc. Roy. Soc. Edinburgh Sect. A \textbf{90} (1981), no. 1-2, 63--70

\bibitem[C82]{C82} B. Carl, On a characterization of operators from $l_q$ into a Banach space of type $p$ with some applications to eigenvalue problems. J. Funct. Anal. \textbf{48} (1982), 394--407

\bibitem[C85]{C85} B. Carl, Inequalities of Bernstein-Jackson type and the degree of compactness of operators in Banach spaces. Ann. Inst. Fourier \textbf{35} (1985), 79--118

\bibitem[C97]{C97} B. Carl, Metric entropy of convex hulls in Hilbert spaces. Bull. London Math. Soc. \textbf{29} (1997), 452--458

\bibitem[CE01]{CE01} B. Carl, D. E. Edmunds, Entropy of $C(X)$-valued operators and diverse applications. J. Inequal. Appl. \textbf{6} (2001), 119--147  

\bibitem[CE03]{CE03} B. Carl, D. E. Edmunds, Gelfand numbers and metric entropy of convex hulls in Hilbert spaces. Studia Math. \textbf{159} (3) (2003), 391--402

\bibitem[CHK88]{CHK88} B. Carl, S. Heinrich, T. Kühn, $s$-numbers of integral operators with Hölder continuous kernels over metric compacta. J. Funct. Anal. \textbf{81} (1988), 54--73

\bibitem[CHP11]{CHP11} B. Carl, A. Hinrichs, A. Pajor, Gelfand numbers and metric entropy of convex hulls in Hilbert spaces. to appear in Positivity

\bibitem[CKP99]{CKP99} B. Carl, I. Kyrezi, A. Pajor, Metric entropy of convex hulls in Banach spaces. J. London Math. Soc. \textbf{60} (2) (1999), 871--896

\bibitem[CP88]{CP88} B. Carl, A. Pajor, Gelfand numbers of operators with values in Hilbert spaces. Invent. Math. \textbf{94} (1988), 459--504

\bibitem[CS90]{CS90} B. Carl, I. Stephani, Entropy, Compactness and Approximation of Operators. Cambridge Univ. Press, Cambridge, 1990

\bibitem[CoKü88]{CoKü88} F. Cobos, T. Kühn, Entropy and eigenvalues of weakly singular integral operators. Integral Equations Operator Theory \textbf{11} (1) (1988), 64--86

\bibitem[CoKü90]{CoKü90} F. Cobos, T. Kühn, Eigenvalues of weakly singular integral operators. J. London Math. Soc. (2) \textbf{41} (2) (1990), 323-335 

\bibitem[CrSt02]{CrSt02} J. Creutzig, I. Steinwart, Metric entropy of convex hulls in type $p$ spaces - the crtitical case. Proc. Amer. Math. Soc. \textbf{130} (2002), 733--743

\bibitem[DJT95]{DJT95} J. Diestel, H. Jarchow, A. Tonge, Absolutely Summing Operators. Cambridge Univ. Press, Cambridge, 1995 

\bibitem[Do93]{Do93} M.R. Dostanic, Asymptotic behavior of the singular values of fractional integral operators. J. Math. Anal. Appl. \textbf{175} (1993), 380--391 

\bibitem[Do95]{Do95} M.R. Dostanic, An estimation of singular values of convolution operators. Proc. Amer. Math. Soc. \textbf{123} (1995), 1399--1409

\bibitem[DM97]{DM97} M.R. Dostanic, D.Z. Milinkovic, Asymptotic behavior of singular values of certain integral operators. Publ. Inst. Math. (Beograd) (N. S.) \textbf{62} (1997), 83--98 

\bibitem[D67]{D67} R.M. Dudley, The sizes of compact subsets of Hilbert space and continuity of Gaussian processes. J. Funct. Anal. \textbf{1} (1967), 290--330

\bibitem[D73]{D73} R.M. Dudley, Sample functions of the Gaussian process. Ann. Probability \textbf{1} (1973), 66--103

\bibitem[D87]{D87} R.M. Dudley, Universal Donsker classes and metric entropy, Ann. Probability \textbf{15} (1987), 1306--1326

\bibitem[ET96]{ET96} D.E. Edmunds, H. Triebel, Function Spaces, Entropy Numbers and Differential Operators. Cambridge Univ. Press, Cambridge, 1996

\bibitem[FM86]{FM86} V. Faber, G.M. Wing, Singular values of fractional integral operators: a unification of theorems of Hille, Tamarkin, and Chang.  J. Math. Anal. Appl. \textbf{120} (2) (1986), 745--760 

\bibitem[G01]{G01} F. Gao, Metric entropy of convex hulls. Israel J. Math. \textbf{123} (2001), 359--364

\bibitem[G04]{G04} F. Gao, Entropy of absolute convex hulls in Hilbert spaces. Bull. London Math. Soc. \textbf{36} (2004), 460--468

\bibitem[G12]{G12} F. Gao, Optimality of CKP-inequality in the critical case, to appear in Proc. Amer. Math. Soc.

\bibitem[GG84]{GG84} A. Garnaev, E. Gluskin, On diameters of the Euclidean sphere. Dokl. Akad. Nauk. USSR \textbf{277} (1984), 1048--1052, Engl. transl.: Soviet Math. Dokl. \textbf{30} (1984), 200--204

\bibitem[Go88]{Go88} Y. Gordon, On Milman's inequality and random subspaces which escape through a mesh in $R^n$. Geometric aspects of functional analysis (1986-87), 84--106, Lecture Notes in Math. \textbf{1317} (1988), Springer Berlin

\bibitem[GKS87]{GKS87} Y. Gordon, H.König and C. Schütt, Geometric and probabilistic estimates for entropy and approximation numbers of operators. J. Approx. Theory \textbf{49} (1987), 219--239

\bibitem[HLP88]{HLP88} G.H. Hardy, J.E. Littlewood, G. Pólya, Inequalities. Reprint of the 1952 edition. Cambridge Mathematical Library. Cambridge University Press, Cambridge, 1988

\bibitem[Ka77]{Ka77} B.S. Kashin, Diameters of some finite-dimensional sets and classes of smooth functions, Math. USSR-Izv. \textbf{11} (1977), 317--333

\bibitem[Kl12a]{Kl12a} O. Kley, Entropy of convex hulls and Kuelbs-Li inequalities. Dissertation (2012), Friedrich-Schiller-University Jena 

\bibitem[Kl12b]{Kl12b} O. Kley, Kuelbs-Li inequalities and metric entropy of convex hulls. J. Theoret. Probab. (2012), DOI: 10.1007/s10959-012-0408-5 

\bibitem[Kö86]{Kö86} H. König, Eigenvalue Distribution of Compact Operators. Birkhäuser, Basel, 1986

\bibitem[Kü05]{Kü05} T. Kühn, Entropy numbers of general diagonal operators. Rev. Mat. Complut. \textbf{18} (2) (2005), 479--491
 
\bibitem[Ky00]{Ky00} I. Kyrezi, On the entropy of the convex hull of finite sets. Proc. Amer. Math. Soc. \textbf{128} (2000), 2393--2403
 
\bibitem[LT91]{LT91} M. Ledoux, M. Talagrand, Probability in Banach Spaces. Springer, Berlin, 1991

\bibitem[Lif10]{Lif10} M.A. Lifshits, Bounds for entropy numbers of some critical operators. Trans. Amer. Math. Soc. \textbf{364} (4) (2012), 1797--1813

\bibitem[Lin04]{Lin04} W. Linde, Kolmogorov numbers of Riemann-Liouville operators over small sets and applications to Gaussian processes. J. Approx. Theory \textbf{128} (2004), 207--233

\bibitem[Lin08]{Lin08} W. Linde, Non-determinism of linear operators and lower entropy estimates, J. Fourier Anal. Appl. \textbf{14} (2008), 568--587

\bibitem[LL99]{LL99} W.V. Li, W. Linde, Approximation, metric entropy and small ball estimates for Gaussian measures. Ann. Probability \textbf{27} (1999), 1556--1578

\bibitem[LL00]{LL00} W.V. Li, W. Linde, Metric entropy of convex hulls in Hilbert spaces. Studia Math. \textbf{139} (1) (2000), 29--45 

\bibitem[LP74]{LP74} W. Linde, A. Pietsch, Mappings of Gaussian measures of cylindrical sets in Banach spaces. Teor. Verojatnost. i Primenen. \textbf{19} (1974), 472--487 

\bibitem[LS06]{LS06} E.N. Lomakina, V.D. Stepanov, Asymptotic estimates for the approximation and entropy numbers of the one-weight Riemann-Liouville operator. (Russian summary) Mat. Tr. \textbf{9} (1) (2006), 52--100

\bibitem[M01]{M01} A. Meskhi, Asymptotic behavior of singular and entropy numbers for some Riemann-Liouville type operators. Georgian Math. J. \textbf{8} (2001), 323--332

\bibitem[PT86]{PT86} A. Pajor, N. Tomczak-Jaegermann, Subspaces of small codimension of finite-dimensional Banach spaces. Proc. Amer. Math. Soc. \textbf{97} (1986), 637--642

\bibitem[P74]{P74} A. Pietsch, $s$-numbers of operators in Banach spaces. Studia Math. \textbf{51} (1974), 201--223

\bibitem[P78]{P78} A. Pietsch, Operator Ideals. VEB Deutscher Verlag der Wissenschaften, Berlin, 1978

\bibitem[P87]{P87} A. Pietsch, Eigenvalues and $s$-Numbers. Cambridge Studies in Advanced Mathematics 13, Cambridge Univ. Press, Cambridge, 1987

\bibitem[P07]{P07} A. Pietsch, History of Banach Spaces and Linear Operators. Birkhäuser Boston, Boston, 2007

\bibitem[Pin85]{Pin85} A. Pinkus, $N$-Widths in Approximation Theory, Springer-Verlag, Berlin, 1985

\bibitem[Pi73a]{Pi73a} G. Pisier, ``Type'' des espaces normés. École Polytech. Palaiseau, Séminaire Maurey-Schwartz 1973/74, Exp. III  

\bibitem[Pi73b]{Pi73b} G. Pisier, Sur les espaces qui ne contiennent pas de $l_1^n$ uniformément. École Polytech. Palaiseau, Séminaire Maurey-Schwartz 1973/74, Exp. VII  

\bibitem[Pi81]{Pi81} G. Pisier, Remarques sur un résultat non publié de B. Maurey. Séminaire d'Analyse Fonctionnelle, École Polytechnique, Palaiseau, Exposé 5, 1980--1981

\bibitem[Pi89]{Pi89} G. Pisier, The Volume of Convex Bodies and Banach Space Geometry. Cambridge Univ. Press, Cambridge, 1989

\bibitem[RS96]{RS96} C. Richter, I. Stephani, Entropy and the approximation of bounded functions and operators. Arch. Math. (Basel) \textbf{67} (6) (1996), 478--492

\bibitem[Sch84]{Sch84} C. Schütt, Entropy Numbers of Diagonal Operators between Symmetric Banach Spaces. J. Approx. Theory \textbf{40} (1984), 121--128

\bibitem[St99]{St99} I. Steinwart, Entropy of $C(K)$-valued operators and some applications. Dissertation (1999), Friedrich-Schiller-University Jena 

\bibitem[St00]{St00} I. Steinwart, Entropy of $C(K)$-valued operators. J. Approx. Theory \textbf{103} (2000), 302--328

\bibitem[St04]{St04} I. Steinwart, Entropy of convex hulls - some Lorentz norm results. J. Approx. Theory \textbf{128} (2004), 42--52

\bibitem[Ta93]{Ta93} M. Talagrand, New Gaussian estimates for enlarged balls. Geom. Funct. Anal. \textbf{3} (1993), 502--526 

\bibitem[To87]{To87} N. Tomczak-Jaegermann, Dualité des nombres d'entropie pour des opérateurs à valeurs dans un espace de Hilbert. C. R. Acad. Sci. Paris Sér. I Math. \textbf{305} (7) (1987), 299--301

\bibitem[Tr70]{Tr70} H. Triebel, Interpolationseigenschaften von Entropie- und Durchmesseridealen kompakter Operatoren. Studia Math. \textbf{34} (1970), 89--107

\end{thebibliography}
\end{document}